\documentclass[12pt]{article}
\usepackage{amsmath,amssymb}
\usepackage{graphicx}
\begin{document}

\title{Computing Nearly Singular Solutions Using Pseudo-Spectral Methods}
\author{Thomas Y. Hou\thanks{Applied and Comput. Math, 217-50, Caltech, Pasadena, 
CA 91125. Email: hou@acm.caltech.edu, and LSEC, 
Academy of Mathematics and Systems Sciences, Chinese Academy of Sciences,
Beijing 100080, China.}
\and Ruo Li\thanks{LMAM\&School of Mathematical Sciences, Peking
  University, Beijing 100871, China.
Email: rli@acm.caltech.edu.}}

\maketitle

\begin{abstract}
In this paper, we investigate the performance of pseudo-spectral
methods in computing nearly singular solutions of fluid dynamics
equations. We consider two different ways of removing the
aliasing errors in a pseudo-spectral method. The first one
is the traditional 2/3 dealiasing rule. The second one 
is a high (36th) order Fourier smoothing which keeps a significant
portion of the Fourier modes beyond the 2/3 cut-off point in the
Fourier spectrum for the 2/3 dealiasing method. Both the 1D Burgers 
equation and the 3D incompressible Euler equations are considered. We 
demonstrate that the pseudo-spectral method with the high order 
Fourier smoothing gives a much better performance than the
pseudo-spectral method with the 2/3 dealiasing rule. Moreover, 
we show that the high 
order Fourier smoothing method captures about $12 \sim 15\%$ 
more effective Fourier modes in each dimension than the 2/3 
dealiasing method. For the 3D Euler equations, the gain in the
effective Fourier codes for the high order Fourier smoothing
method can be as large as 20\% over the 2/3 dealiasing method. 
Another interesting observation is that the error produced by 
the high order Fourier smoothing method is highly localized near 
the region where the solution is most singular, while the 2/3 
dealiasing method tends to produce oscillations in the entire 
domain. The high order Fourier
smoothing method is also found be very stable dynamically.
No high frequency instability has been observed.
\end{abstract}

\section{Introduction}

Pseudo-spectral methods have been one of the most commonly used
numerical methods in solving nonlinear partial differential
equations with periodic boundary conditions. Pseudo-spectral
methods have the advantage of computing a nonlinear convection
term very efficiently using the Fast Fourier Transform. On the
other hand, the discrete Fourier transform of a periodic
function introduces the so-called aliasing error
\cite{GO77,CHQZ88,CHQZ06, Boyd00},
which is partially due to the artificial periodicity of the 
discrete Fourier coefficient as a function of the wave number. 
The aliasing error pollutes
the accuracy of the high frequency modes, especially those
last 1/3 of the high frequency modes.  Without using any 
dealiasing or Fourier smoothing, the pseudo-spectral method may 
suffer from some mild numerical instability \cite{GHT94}.  
One of the most commonly used dealiasing methods is the 
so-called 2/3 dealiasing rule, in which one sets to zero
the last 1/3 of the high frequency modes and keeps 
the first 2/3 of the Fourier modes unchanged.
Another way to control the aliasing errors is to apply
a smooth cut-off function or Fourier smoothing to the 
Fourier coefficients. However, many existing Fourier 
smoothing methods damp the last 1/3 of the high frequency 
modes just like the 2/3 dealiasing method.

In this paper, we investigate the performance of pseudo-spectral
methods using the 2/3 dealiasing rule and a high order Fourier 
smoothing. In the Fourier smoothing method, we use a 36th order 
Fourier smoothing function which keeps a significant portion of 
the Fourier modes beyond the 2/3 cut-off point in the Fourier
spectrum for the 2/3 dealiasing 
rule. We apply these two methods to compute nearly singular 
solutions in fluid flows. Both the 1D Burgers equation
and the 3D incompressible Euler equations will be considered.
The advantage of using the Burgers equation is that it shares some
essential difficulties as other fluid dynamics equations, and yet we 
have a semi-analytic formulation for its solution. By using the 
Newton iterative method, we can obtain an approximate solution 
to the exact solution up to 13 digits of accuracy. Moreover, 
we know exactly when a shock singularity will form in time. This 
enables us to perform a careful convergence study in both the 
physical space and the spectral space very close to the singularity 
time. 

We first perform a careful convergence study of the two pseudo-spectral 
methods in both physical and spectral spaces for the 1D Burgers 
equation. Our extensive numerical results demonstrate that the 
pseudo-spectral method with the high order Fourier smoothing
(the Fourier smoothing method for short) gives a much better performance 
than the pseudo-spectral method with the 2/3 dealiasing rule
(the 2/3 dealiasing method for short). In particular, we show that 
the unfiltered high frequency coefficients in the Fourier smoothing method 
approximate accurately the corresponding exact Fourier coefficients. 
More precisely, we demonstrate that the Fourier 
smoothing method captures about $12 \sim 15\%$ more effective Fourier 
modes than the 2/3 dealiasing method in each dimension. The gain is 
even higher for the 3D Euler equations since the number of effective 
modes in the Fourier smoothing method is higher in three 
dimensions. Thus the Fourier smoothing method gives a 
more accurate approximation than the 2/3 dealiasing method. We 
will illustrate this improved accuracy by studying the errors in
$L^\infty$-norm and $L^1$-norm as a function of time, and by studying
the spatial distribution of the pointwise error and the convergence of 
the Fourier spectrum at a sequence of times very close to the singularity
time. Another interesting observation is that the error produced by 
the Fourier smoothing method is highly localized near the 
region where the solution is most singular and decays exponentially fast 
with respect to the distance from the singularity point. The error in 
the smooth region is several orders of magnitude smaller than that in 
the singular region. On the other hand, the 
2/3 dealiasing method produces noticeable oscillations in the entire 
domain as we approach the singularity time. This is to some extent
due to the Gibbs phenomenon and the loss of the $L^2$ energy associated
with the solution.
Moreover, our computational results show that in the smooth region
the error produced by the Fourier smoothing method is several orders
of magnitude smaller than that produced by the 2/3 dealiasing method.
This is an important advantage of the Fourier smoothing 
method over the 2/3 dealiasing method.   

Next, we apply the two pseudo-spectral methods to solve the 
nearly singular solution of the 3D incompressible Euler equations.
We would like to see if the comparison we make regarding the convergence
properties of the two pseudo-spectral methods for the 1D Burgers equation
is still valid for the more challenging 3D incompressible Euler equations. 
In order to make our comparison meaningful, we choose a smooth initial 
condition which could potentially develop a finite time singularity.
There have been many computational efforts in searching for finite time
singularities of the 3D Euler equations, see e.g.
\cite{Chorin82,PS90,KH89,GS91,SMO93,Kerr93,Caf93,BP94,FZG95,Pelz98,GMG98,Kerr05}.
One of the frequently cited numerical evidences for a finite time blowup
of the Euler equations is the two slightly perturbed anti-parallel vortex tubes
initial data studied by Kerr \cite{Kerr93,Kerr05}. In Kerr's computations, 
a pseudo-spectral discretization with the 2/3 dealiasing rule was used 
in the $x$ and $y$ directions while a Chebyshev polynomial discretization 
was used along the $z$ direction. His best space resolution was of the 
order $512\times 256\times 192$. In \cite{Kerr93,Kerr05}, Kerr reported 
that the maximum vorticity blows up like $O((T-t)^{-1})$ and the velocity 
field blows up like $O((T-t)^{-1/2})$. The alleged singularity time $T$ 
is equal to 18.7 while his computations beyond $t=17$ were not considered 
as the primary evidence for a singularity since they were polluted by 
noises \cite{Kerr93}.

We perform a careful convergence study of the two pseudo-spectral methods
using Kerr's initial condition with a sequence of resolutions up to $T=19$, 
beyond the singularity time alleged in \cite{Kerr93,Kerr05}. The largest 
space resolution we use is $1536\times 1024 \times 3072$. Convergence in 
both physical and spectral space has been observed for the two
pseudo-spectral methods. Both numerical methods converge to the same 
solution under mesh refinement. Our computational study also demonstrates
that the Fourier smoothing method offers better computational accuracy 
than the 2/3 dealiasing method. For a given resolution, the Fourier 
smoothing method captures about 20\% more effective Fourier modes 
than the 2/3 dealiasing method does. We also find that the 2/3 
dealiasing method produces some oscillations near the 2/3 
cut-off point of the spectrum. This abrupt cut-off of the Fourier 
spectrum generates noticeable oscillations in the vorticity contours 
at later times. Even using a relative high resolution 
$1024\times 786 \times 2048$, we find that the vorticity contours 
obtained by the 2/3 dealiasing method still suffer relatively large 
oscillations in the late stage of the computations, which are to
some extent caused by the Gibbs phenomenon and the loss of enstrophy 
due to the abrupt cut-off of the high frequency modes.
On the other hand, the vorticity contours 
obtained by the Fourier smoothing method remains smooth throughout the 
computations. Our spectral computations using both the 2/3 dealiasing 
rule and the high order Fourier smoothing confirm the finding 
reported in \cite{HL06}, i.e. the maximum vorticity does not grow 
faster than double exponential in time and the velocity field remains 
bounded up to $T=19$.

We would like to emphasize that the Fourier smoothing method is very 
stable and robust in all our computational experiments. We do not observe 
any high frequency instability in our computations for both the 1D Burgers 
equation and the 3D incompressible Euler equations. The resolution study 
that we conduct is completely based on the consideration of accuracy, not 
by the consideration of stability. 

We would like to mention that Fourier smoothing has been also used 
effectively to approximate discontinuous solutions of linear hyperbolic
equations, see e.g. \cite{MMO77,ML78,AGT86}. To compute discontinuous
solutions for nonlinear conservation laws, the spectral viscosity method 
has been introduced and analyzed, see \cite{Tad89,MT89} and the review 
article \cite{Tad93}. Like Fourier smoothing, the purpose of introducing
the spectral viscosity is to localize the Gibbs oscillations, maintaining
stability without loss of spectral accuracy.

The remaining of the paper is organized as follows. In Section 2, we 
present a careful convergence study of the two pseudo-spectral methods 
for the 1D Burgers equation. In Section 3, we present a similar convergence
study for the incompressible 3D Euler equations using Kerr's initial 
data. Some concluding remarks are made in Section 4.

\section{Convergence study of the two pseudo-spectral methods for the 1D Burgers equation}

In this section, we perform a careful convergence study of the 
two pseudo-spectral methods for the 1D Burgers equation. The 1D Burgers
equation shares some of the essential difficulties in many fluid
dynamic equations. In particular, it has the same type of quadratic 
nonlinear convection term as other fluid dynamics equations.
It is well known that the 1D Burgers equation can form a shock 
discontinuity in a finite time \cite{LeVeque92}. The advantage 
of using the 1D Burgers equation as a prototype is that we have 
a semi-analytical solution formulation for the 1D Burgers equation. 
This allows us to use the Newton iterative method to obtain a very accurate 
approximation (up to 13 digits of accuracy) to the exact solution of 
the 1D Burgers equation arbitrarily close to the singularity time. This 
provides a solid foundation in our convergence study of the two
spectral methods.

We consider the inviscid 1D Burgers equation 
\begin{equation}\label{eq.burgers}
u_t + \left(\frac{u^2}{2}\right)_x = 0, \quad -\pi \leq x \leq \pi,
\end{equation}
with an initial condition given by
\[ 
u|_{t=0} = u_0(x). 
\]
We impose a periodic boundary condition over $[-\pi, \pi]$. 
By the method of characteristics,
it is easy to show that the solution of the 1D Burgers equation is given by
\begin{equation}
\label{burgers-sol}
u(x,t)=u_0(x-t u(x,t)) .
\end{equation}
The above implicit formulation defines a unique solution for 
$u(x,t)$ up to the time when the first shock singularity 
develops. After the shock singularity develops, equation (\ref{burgers-sol})
gives a multi-valued solution. An entropy condition is required to
select a unique physical solution beyond the shock singularity \cite{LeVeque92}. 

We now use a standard pseudo-spectral method to
approximate the solution. Let $N$ be an integer, and let 
$h=\pi/N$. We denote by $x_j = j h$ ($j=-N,...,N$) the discrete 
mesh over the interval $[-\pi,\pi]$. To describe the pseudo-spectral
methods, we recall that the discrete Fourier transform of a periodic
function $u(x)$ with period $2 \pi$ is defined by
\[ 
\hat{u}_k = \frac{1}{2N}\sum_{j=-N+1}^{N} u(x_j) e^{-i k x_j} \;.
\]
The inversion formula reads
\[
u(x_j) = \sum_{k=-N+1}^{N} \hat{u}_k e^{i k x_j} \; .
\]
We note that $\hat{u}_k $ is periodic in $k$ with period $2N$. This
is an artifact of the discrete Fourier transform, and the source of
the aliasing error. To remove the aliasing error, one usually applies
some kind of dealiasing filtering when we compute the discrete 
derivative. Let $\rho(k/N)$ be a cut-off function in the spectrum
space. A discrete derivative operator may be expressed in the
Fourier transform as 
\begin{equation}
\label{spec_deriv}
\widehat{(D_h u)}_k = ik \rho (k/N) \widehat{u}_k , \quad k=-N+1, ...,N.
\end{equation}
Both the 2/3 dealiasing rule and the Fourier
smoothing method can be described by a specific choice of the
high frequency cut-off function, $\rho$ (also known as Fourier
filter). For the 2/3 dealiasing rule, the cut-off function is 
chosen to be 
\begin{equation}
\label{cutoff-23rd}
\rho (k/N) = \left \{ \begin{array}{ll} 1, & \mbox{if $ |k/N| \leq 2/3$},\\
                                      0, & \mbox{if $ |k/N| > 2/3$.} 
                    \end{array}
            \right .
\end{equation}
In our computations, in order to obtain an alias-free computation on
a grid of $M$ points for a quadratic nonlinear equation, we apply the 
above filter to the high wavenumbers so as to retain only $(2/3)M$ 
unfiltered wavenumbers before making the coefficient-to-grid Fast 
Fourier Transform. This dealiasing procedure is alternatively known
as the the 3/2 dealiasing rule because to obtain $M$ unfiltered 
wavenumbers, one must compute nonlinear products in physical space on 
a grid of $(3/2)M$ points, see page 229 of \cite{Boyd00} for more 
discussions.

For the Fourier smoothing method, we choose $\rho$ as follows:
\begin{equation}
\label{cutoff-sm}
\rho (k/N) =  e^{-\alpha (|k|/N)^m},
\end{equation}
with $\alpha = 36$ and $m=36$. In our implementation, both filters are
applied on the numerical solution at every time step. For the $2/3$ 
dealiasing rule, the Fourier modes with wavenumbers $|k| \ge 2/3N$ 
are always set to zero. Thus there is no aliasing error being 
introduced in our approximation of the nonlinear convection term.

The Fourier smoothing method we choose is based on three
considerations.  The first one is that the aliasing instability is 
introduced by the highest frequency Fourier modes. As demonstrated 
in \cite{GHT94}, as long as one can damp out a small portion of 
the highest frequency Fourier modes, the mild instability
caused by the aliasing error can be under control. The second 
observation is that the magnitude of the Fourier coefficient is 
decreasing with respect to the wave number $|k|$ for a function 
that has certain degree of regularity. Typically, we have
$|\widehat{u}_k| \leq C/(1+|k|^m)$ if the $m$th derivative of a 
function $u$ is bounded in $L^1$.
Thus the high frequency Fourier modes have a relatively smaller
contribution to the overall solution than the low to intermediate
frequency modes. The third observation is that one should not 
cut off high frequency Fourier modes abruptly to avoid the 
Gibbs phenomenon and the loss of the $L^2$ energy associated with the 
solution. This is especially important when we compute a nearly 
singular solution whose high frequency Fourier coefficient 
has a very slow decay. 

\begin{figure}
\begin{center}
\includegraphics[width=8cm]{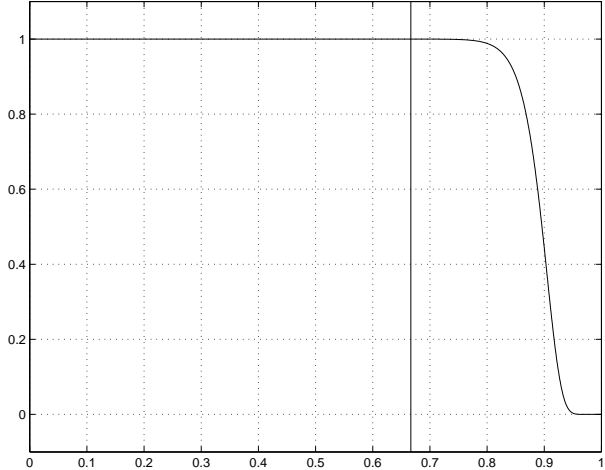}
\end{center}
\caption{The profile of the Fourier smoothing, $\exp(-36 (x)^{36})$,
as a function of $x$. The vertical line corresponds to the cut-off point 
in the Fourier spectrum in the $2/3$ 
dealiasing rule. We can see that using this Fourier smoothing we keep 
about $12 \sim 15\%$ more modes than those using the 2/3 dealiasing rule.
\label{fig.fourier_smoother}}
\end{figure}

Based on the above considerations, we choose a smooth cut-off 
function which 
decays exponentially fast with respect to the high wave number. In 
our cut-off function, we choose the parameters $\alpha=36$ and $m=36$. 
These two parameters are chosen to achieve two objectives: (i) When 
$|k|$ is close to $N$, the cut-off function reaches the machine 
precision, i.e. $10^{-16}$; (ii) The cut-off function remains very 
close to 1 for $|k| < 4N/5$, and decays rapidly and smoothly to zero 
beyond $|k| = 4N/5$. In Figure \ref{fig.fourier_smoother}, we plot
the cut-off function $\rho (x)$ as a function of $x$. The cut-off 
function used by the 2/3rd dealiasing rule is plotted on top of
the cut-off function used by the Fourier smoothing method.
We can see that the Fourier smoothing method keeps about 
$12 \sim 15\%$ more modes than the 2/3 dealiasing method. In
this paper, we will demonstrate by our numerical experiments 
that the extra modes we keep by the Fourier smoothing method give 
an accurate approximation of the correct high frequency Fourier 
modes. 

Next, we will present a convergence study of the two pseudo-spectral methods using
a generic initial condition, $u_0(x) = \sin(x)$. For this initial condition, 
the solution will develop a shock singularity at $t=1$ at $x=\pi$ and $x=-\pi$.
For $t < 1$ but sufficiently close to the singularity time, a sharp layer will 
develop at $x=\pi$ and $x=-\pi$. We will perform a sequence of well-resolved 
computations sufficiently close to the singularity time using the two 
pseudo-spectral methods and study their convergence properties. In our 
computations for both methods, we will use a standard compact three 
step Runge-Kutta scheme for the time integration. In order to compute the 
errors of the two pseudo-spectral methods accurately, we need to solve for the 
``exact solution'' on a computational grid. We do this by solving the 
implicit solution formula (\ref{burgers-sol}) using the Newton iterative 
method up to 13 digits of accuracy.

\begin{figure}
\begin{center}
\includegraphics[width=8cm]{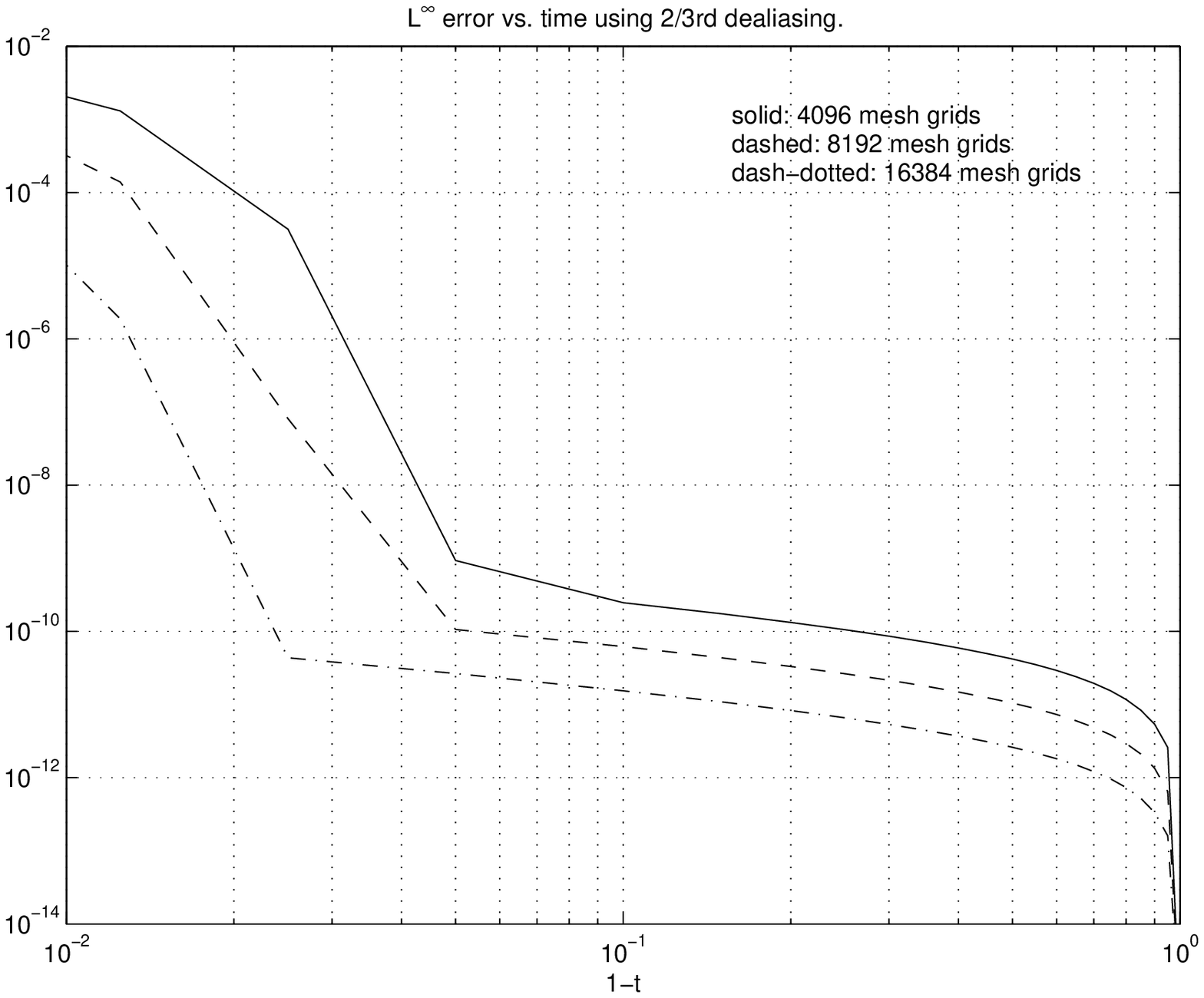}
\includegraphics[width=8cm]{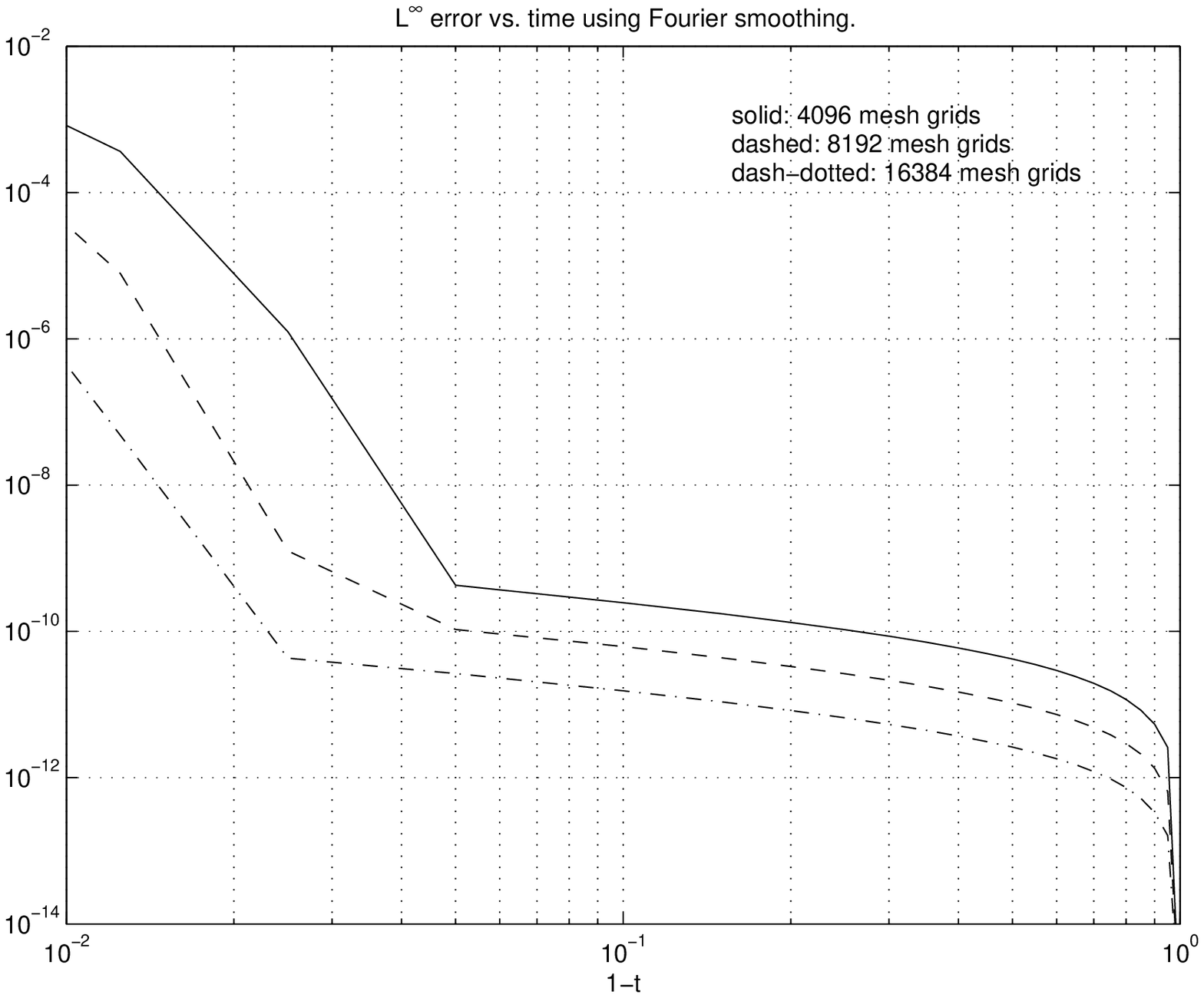}
\end{center}
\caption{The $L^\infty$ errors of the two pseudo-spectral methods as a function
of time using three different resolutions. The plot is in a log-log scale. 
The initial condition is $u_0(x) = \sin(x)$.
The left figure is the result obtained by the 2/3 dealiasing method and 
the right figure is the result obtained by the Fourier smoothing method. 
It can be seen clearly that the $L^\infty$
error obtained by the Fourier smoothing method is smaller than that
obtained by the 2/3 dealiasing method. 
  \label{fig.burgers_l0_error}}
\end{figure}

\begin{figure}
\begin{center}
\includegraphics[width=8cm]{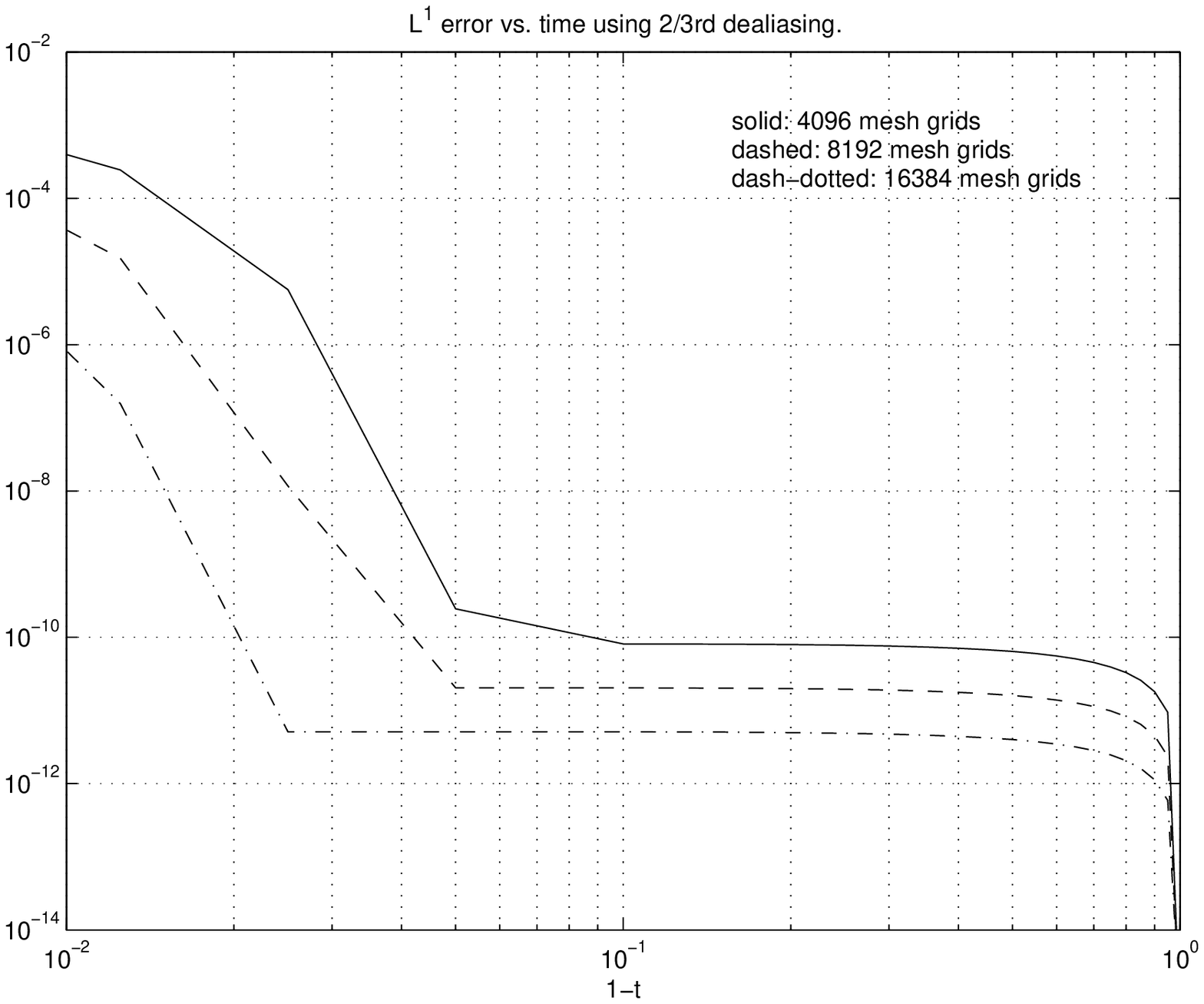}
\includegraphics[width=8cm]{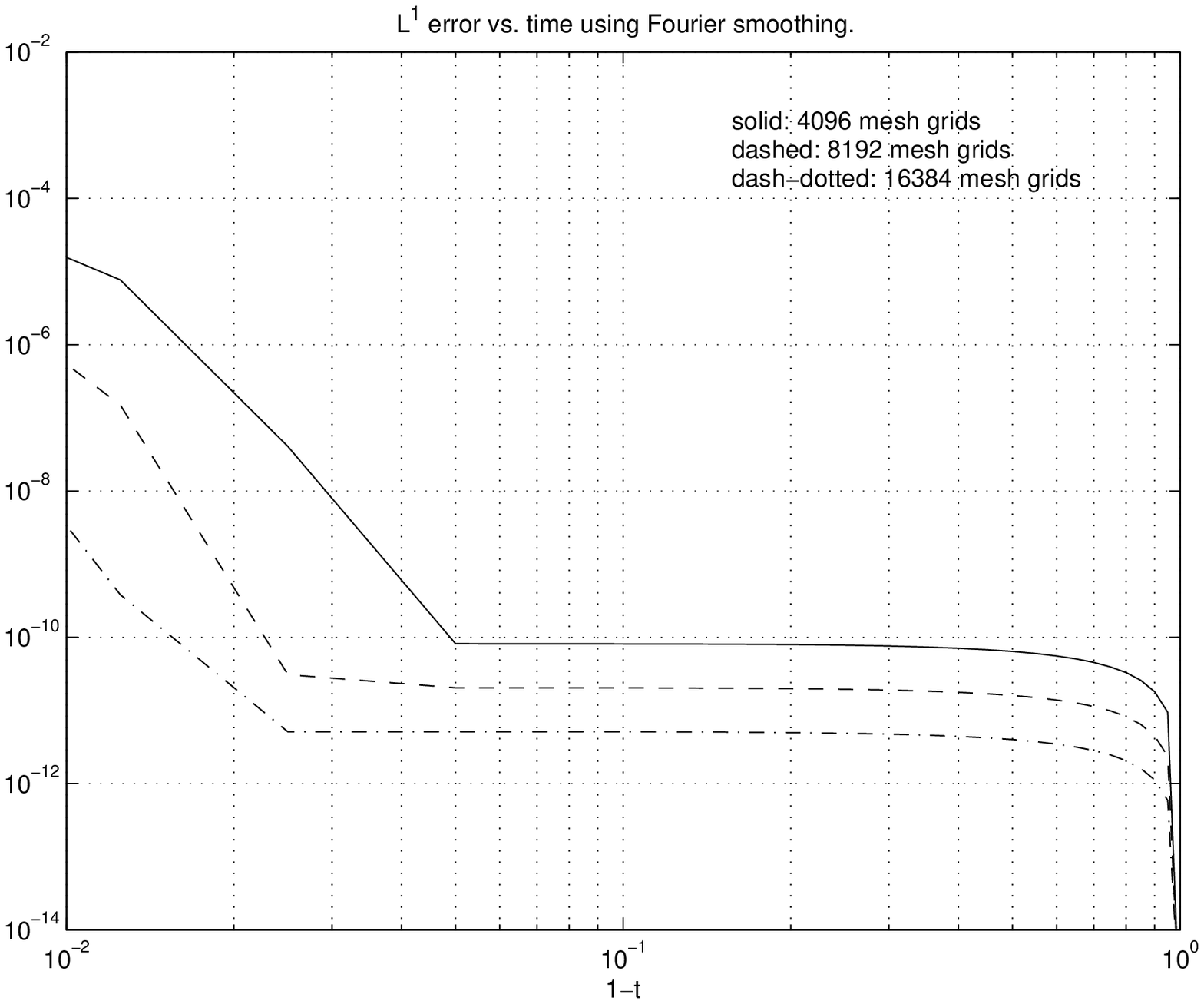}
\end{center}
\caption{The $L^1$ errors of the two pseudo- spectral methods as a function
of time using three different resolutions. The plot is in a log-log scale. 
The initial condition is given by $u_0(x) = \sin(x)$.
  The left figure is the result obtained by the 2/3
  dealiasing method and the right figure is the result obtained by the
  Fourier smoothing method. One can see that the $L^1$ error obtained by 
  the Fourier smoothing method is much smaller than corresponding $L^\infty$ error. 
\label{fig.burgers_l1_error}}
\end{figure}

Our computational studies demonstrate convincingly that the 
Fourier smoothing method gives a much better performance than the
2/3 dealiasing method. In Figure \ref{fig.burgers_l0_error}, we plot the 
$L^\infty$ error of the two pseudo-spectral methods as a function of time 
using three different resolutions. The errors are plotted in a log-log scale. 
The left figure is the result obtained by the 2/3 dealiasing method and the 
right figure is the result obtained by the Fourier smoothing method. We can 
see clearly that the $L^\infty$ error obtained by the Fourier smoothing 
method is smaller than that obtained by the 2/3rd dealiasing method. With 
increasing resolution, the errors obtained by both methods decay rapidly, 
confirming the spectral convergence of both methods. In Figure 
\ref{fig.burgers_l1_error}, we plot the $L^1$ errors
of the two methods as a function of time using three different
resolutions. We can see that the convergence of the Fourier smoothing
method is much faster than the 2/3 dealiasing method. 

\begin{figure}
\begin{center}
\includegraphics[width=8cm]{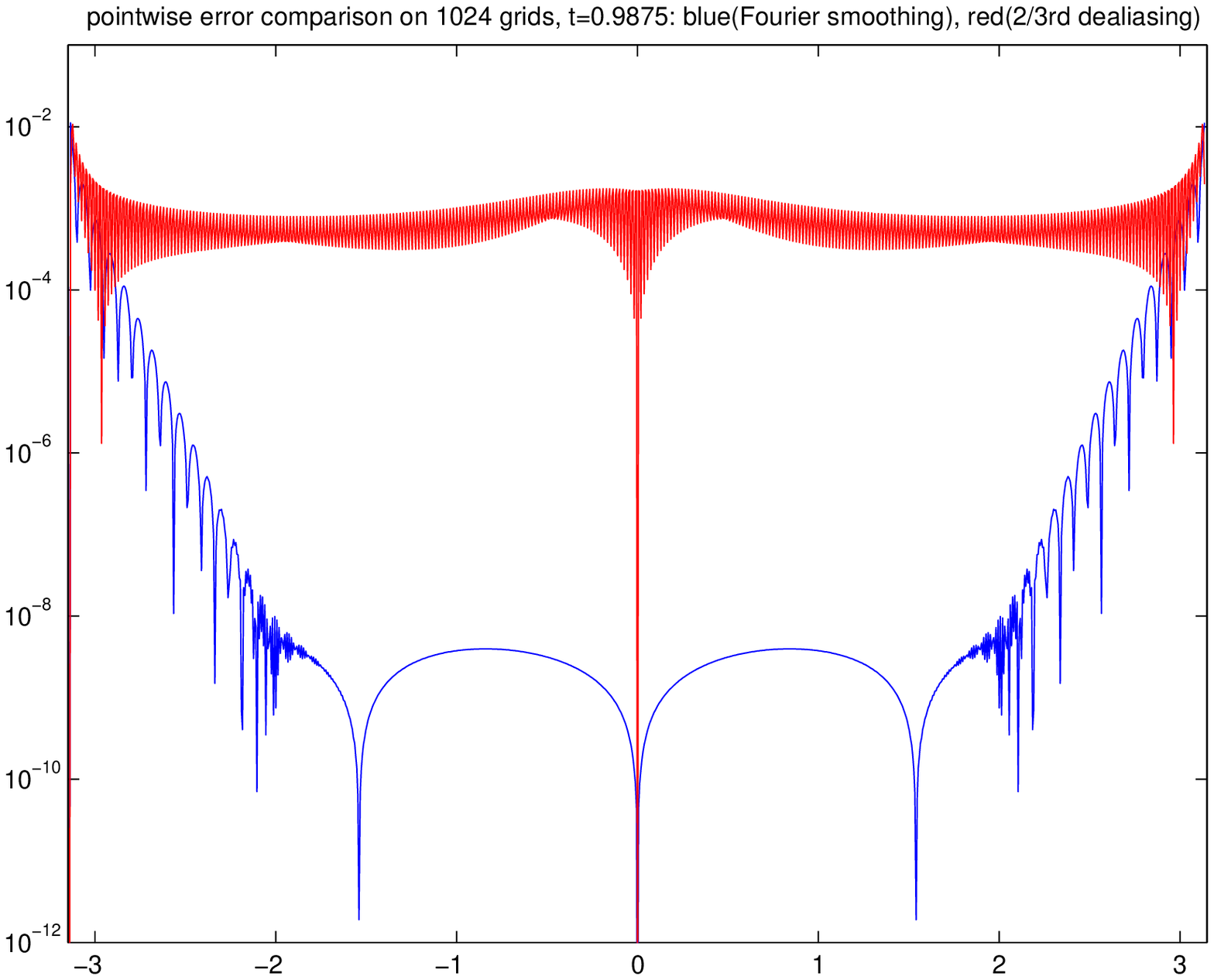}
\includegraphics[width=8cm]{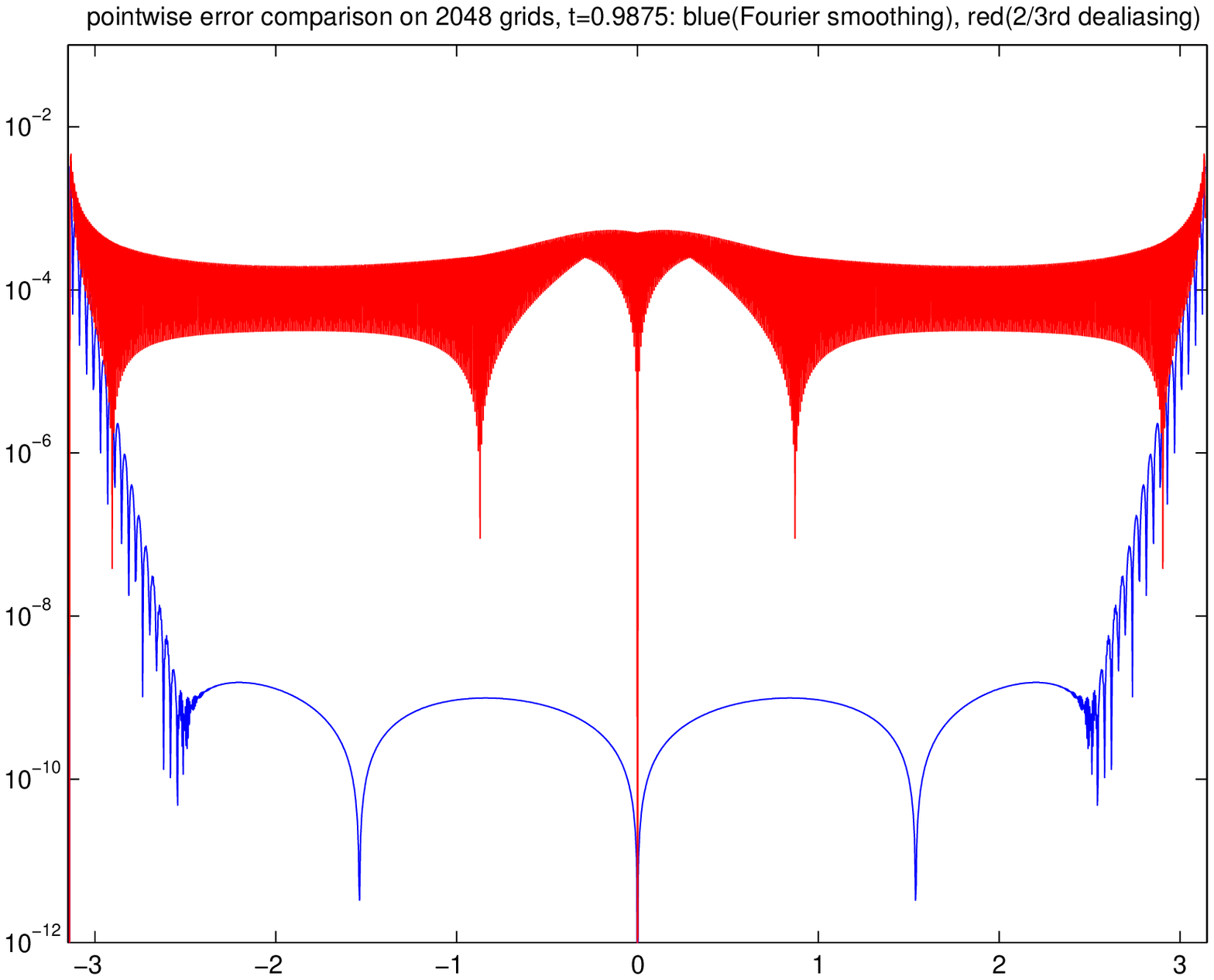}
\end{center}
\caption{The pointwise errors of the two pseudo-spectral methods as a function
of time using three different resolutions. The plot is in a log scale. 
The initial condition is given by $u_0(x) = \sin(x)$.
  The error of the 2/3rd dealiasing method is highly oscillatory and 
  spreads out over the entire domain, while the error of the Fourier 
  smoothing method is highly localized near the location of the shock 
  singularity. 
\label{fig.burgers_pointwise_error}} 
\end{figure}

It is interesting to study the spatial distribution of the pointwise 
errors obtained by the two methods. In Figure 
\ref{fig.burgers_pointwise_error}, we plot the pointwise errors of the 
two methods over one period $[-\pi,\pi]$ at $t=0.985$. In the
computations presented on the left picture, we use resolution $N=1024$, 
while the computations in the right picture corresponds 
to $N=2048$. The errors are plotted in a log scale. We can see that the
error of the 2/3rd dealiasing method, which is colored in red, is highly 
oscillatory and spreads out over the entire domain. This is caused by the
Gibbs phenomenon and the loss of the $L^2$ energy associated with the solution.
On the other hand, the error of the Fourier smoothing method
is highly localized near the location of the shock singularity at $x=\pi$
and $x=-\pi$, and decays exponentially fast with respect to the distance 
from the singularity point. The error in the smooth region is several 
orders of magnitude smaller than that in the singular region. This is 
a very interesting phenomenon. This property makes the Fourier smoothing 
method a better method for computing nearly singular solutions.

\begin{figure}
\begin{center}
\includegraphics[width=8cm]{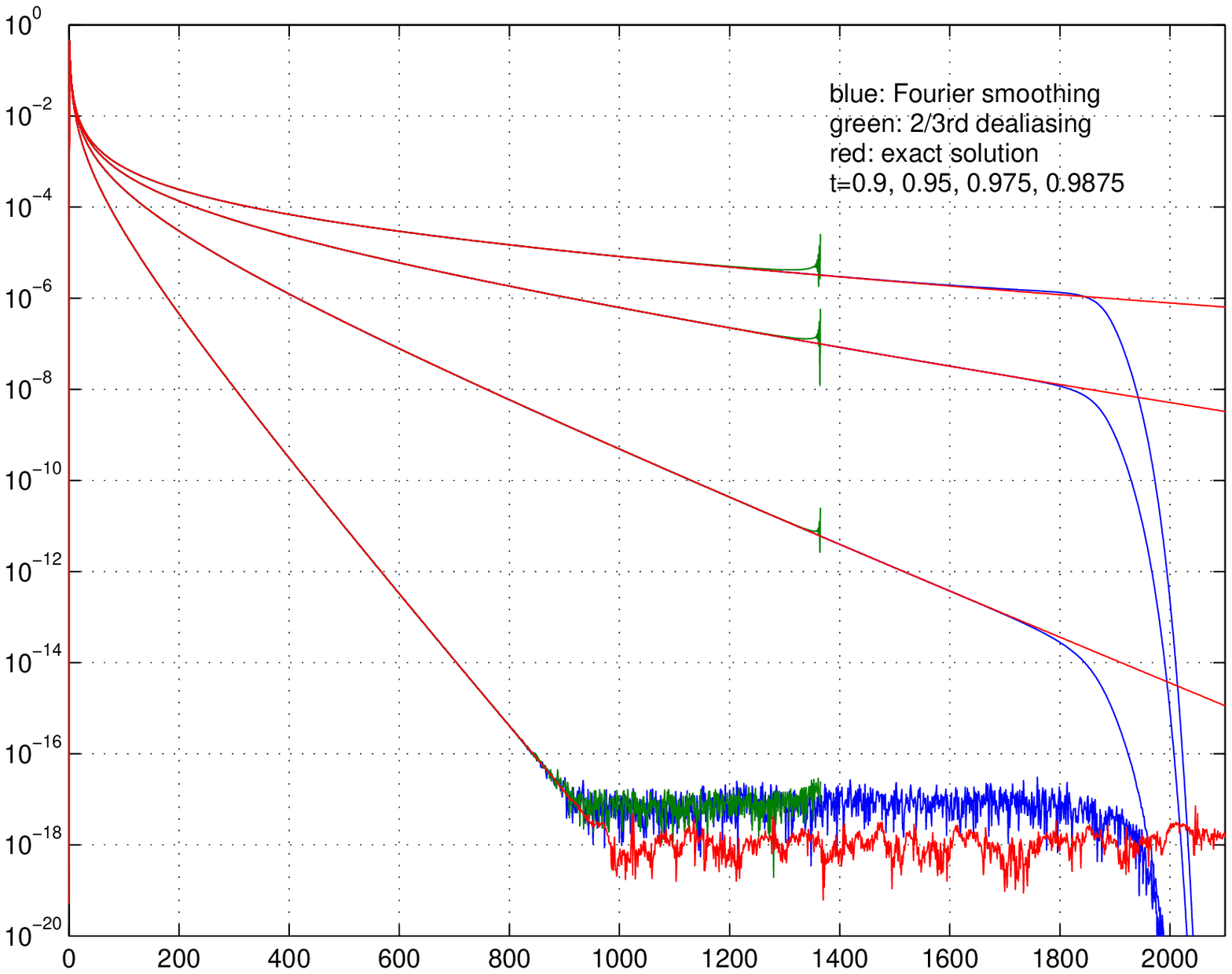}
\includegraphics[width=8cm]{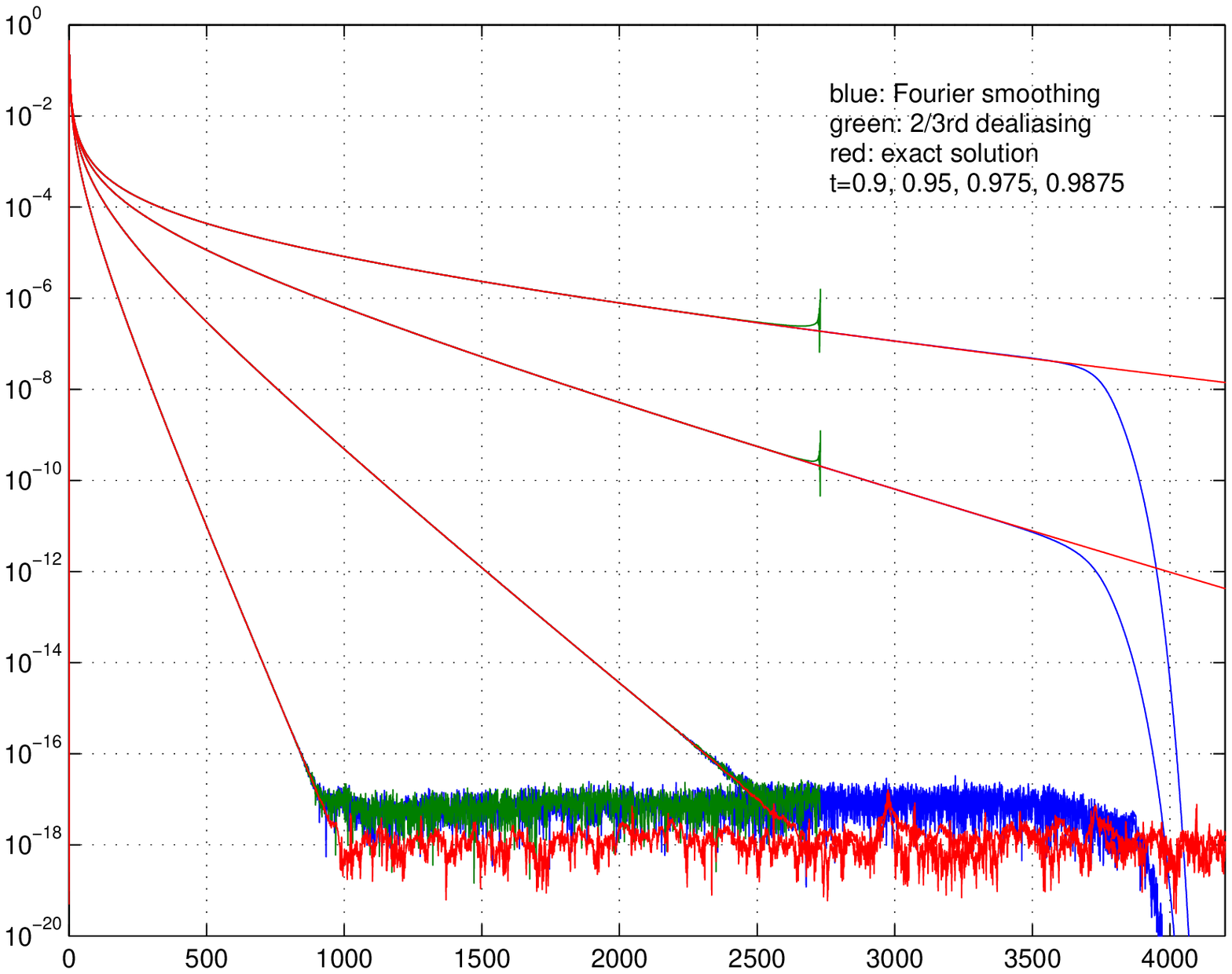}
\end{center}
\caption{Comparison of Fourier spectra of the two methods on different 
resolutions at a sequence of times. The initial condition is given
by $u_0(x) = \sin(x)$.
The left picture corresponds to $N=4096$ and the right picture 
corresponds to $N=8192$.
\label{fig.burgers_spec}}
\end{figure}

To gain further insight of the two pseudo-spectral methods, we study 
the convergence of the two methods in the spectral space. In Figure 
\ref{fig.burgers_spec}, we plot the Fourier spectra of the two
spectral methods at a sequence of times with two different
resolutions $N=4096$ and $N=8192$ respectively. We observe 
that when the solution is relatively smooth and can be resolved
by the computational grid, the spectra of the two methods are
almost indistinguishable. However, when the solution becomes
more singular and cannot be completely resolved by the computational
grid, the two methods give a very different performance. 
For a given resolution, the Fourier smoothing method
keeps about 20\% more Fourier modes than the 2/3 dealiasing method.
When we compare with the ``exact'' spectrum obtained by using
the Newton iterative method, we can see clearly that the extra
Fourier modes that are kept by the Fourier smoothing method give 
indeed an accurate approximation to the correct Fourier modes.
This explains why the Fourier smoothing method offers better
accuracy than the 2/3 dealiasing method. On the other hand,
we observe that the Fourier spectrum of the 2/3 dealiasing method
develops noticeable oscillations near the 2/3 cut-off point of
the Fourier spectrum. This abrupt cut-off in the high frequency
spectrum gives rise to the well-known Gibbs phenomenon and the
loss of the $L^2$ energy, which is 
the main cause for the highly oscillatory and widespread pointwise 
error that we observe in Figure \ref{fig.burgers_pointwise_error}.

\begin{figure}
\begin{center}
\includegraphics[width=8cm]{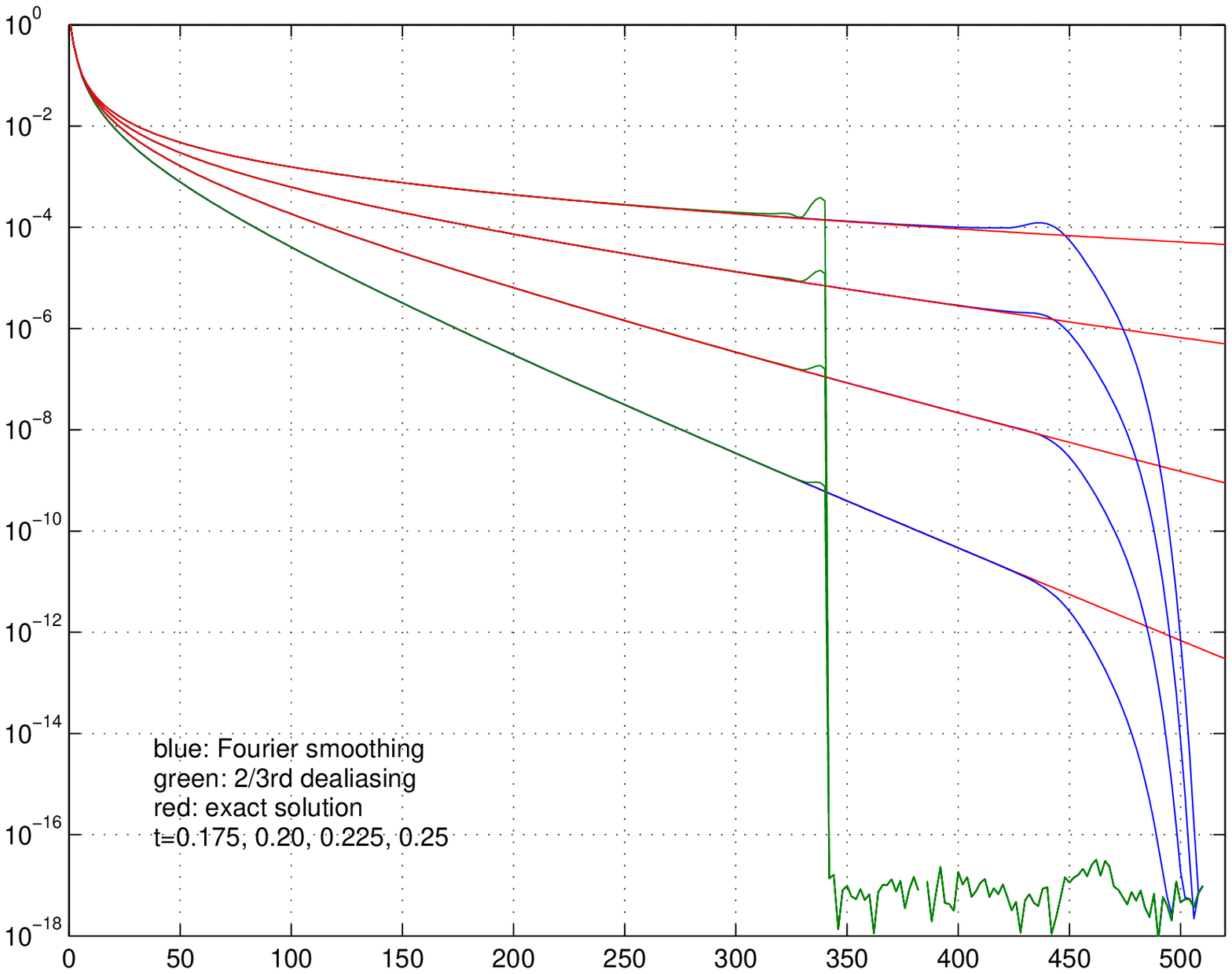}
\includegraphics[width=8cm]{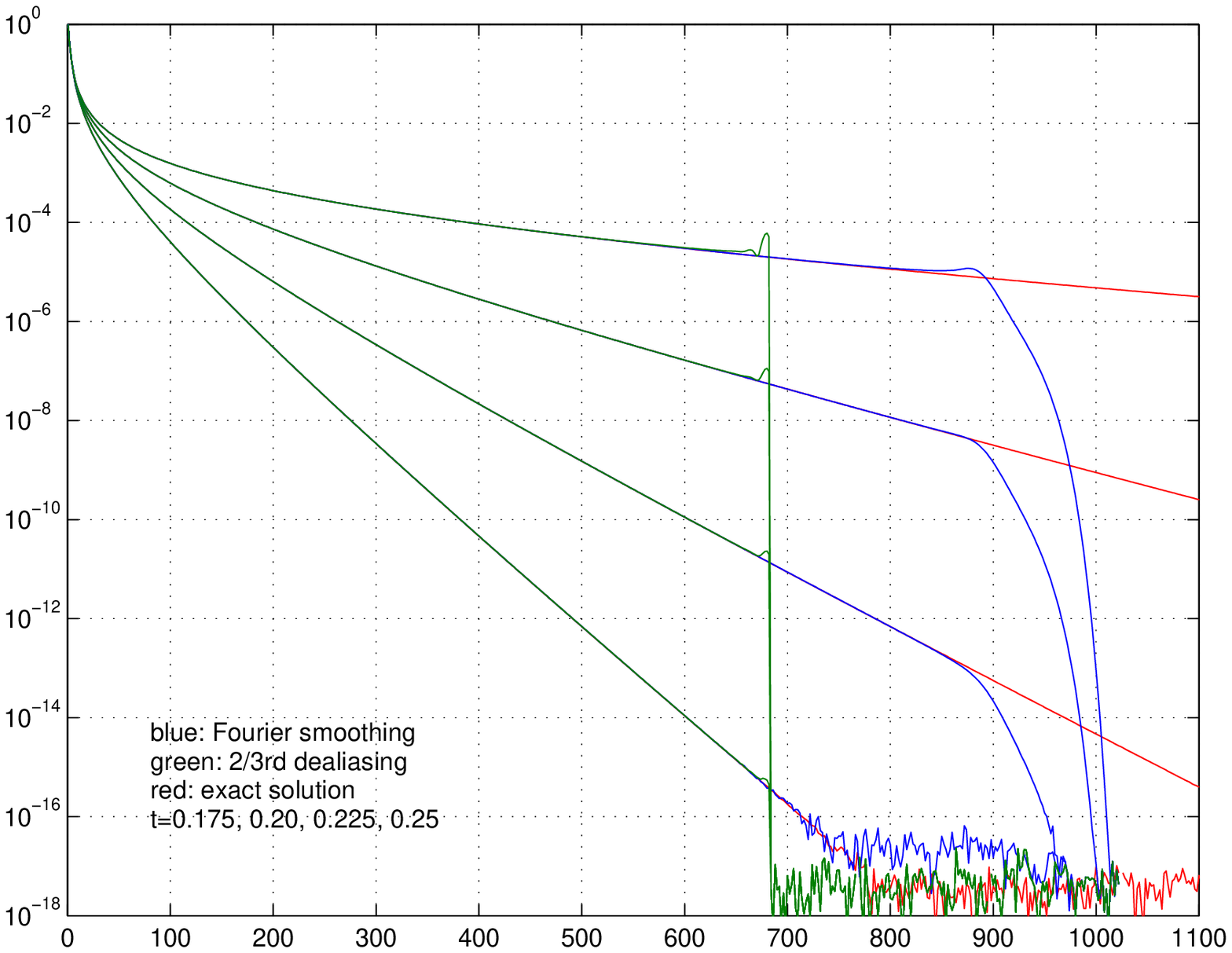}
\end{center}
\caption{Comparison of Fourier spectra of the two methods on different 
resolutions at a sequence of times. The initial condition is given by 
$u_0(x) = (0.1 + \sin^2 x)^{-1/2}$.
The left picture corresponds to $N=1024$ and the right picture 
corresponds to $N=2048$. 
Notice that only even modes are plotted in these figures since the odd
modes are vanished in this example.
\label{fig.ex1_spec}}
\end{figure}

We have also performed similar numerical experiments for several
other initial data. They all give the same qualitative behavior as the
one we have demonstrated above.  In Figure \ref{fig.ex1_spec}, we plot
the Fourier spectra of the two methods at a sequence of times using
a different initial condition: 
$u_0(x) = (0.1 + \sin^2 x )^{-1/2}$. The picture on the left
corresponds to resolution $N=1024$, while the picture on the right
corresponds to resolution $N=2048$. One can see that the convergence
properties of the two methods are essentially the same as those 
presented for the initial condition $u_0(x)=\sin(x)$.

Finally, we would like to point out that the numerical computations
using the Fourier smoothing method have been very stable and robust.
No high frequency instability has been observed throughout our computations.
This indicates that the high order Fourier smoothing we use has effectively
eliminated the mild numerical instability introduced by the 
aliasing error \cite{GHT94}.

%
%

\section{Computing nearly singular solutions of the 3D Euler equations using
pseudo-spectral methods}

In this section, we will apply the two pseudo-spectral methods to solve
the nearly singular solution of the 3D incompressible Euler equations.
The spectral computation of the 3D incompressible Euler equations is
much more challenging due to the nonlocal and nonlinear nature of the
problem and the possible formation of a finite time singularity.
It would be interesting
to find out if the comparison we have made regarding the convergence 
property of the two pseudo-spectral methods for the 1D Burgers equation
is still valid for the 3D Euler equations. To make our comparison useful,
we choose a smooth initial condition which could potentially develop a 
finite time singularity. There have been many computational efforts in 
searching for finite time singularities of the 3D Euler 
equations, see e.g.
\cite{Chorin82,PS90,KH89,GS91,SMO93,Kerr93,Caf93,BP94,FZG95,Pelz98,GMG98,Kerr05}.
Of particular interest is the numerical study of the interaction of two 
perturbed antiparallel vortex tubes by Kerr \cite{Kerr93,Kerr05}, in which 
a finite time blowup of the 3D Euler equations was reported. In this
section, we will perform the comparison of the two pseudo-spectral
methods using Kerr's initial condition.

The 3D incompressible Euler equations in the vorticity stream function
formulation are given as follows (see, e.g., \cite{CM93,MB02}):
\begin{eqnarray}\label{3deuler}
\vec{\omega}_t+(\vec{u}\cdot\nabla) \vec{\omega} & = & \nabla 
\vec{u} \cdot \vec{\omega}, \\
- \bigtriangleup \vec{ \psi} &= & \vec{\omega}, \quad
\vec{u} = \nabla \times \vec{\psi},
\end{eqnarray} 
with initial condition $\vec{\omega}\mid_{t=0} =  \vec{\omega}_{0}$,
where $\vec{u}$ is velocity, $\vec{\omega}$ is vorticity, and 
$\vec{\psi}$ is stream function. Vorticity is related to velocity 
by $\vec{\omega} = \nabla \times \vec{u}$. The incompressibility
implies that 
\[
\nabla \cdot \vec{u} = \nabla \cdot \vec{\omega} = \nabla \cdot \vec{\psi} = 0.
\]
We consider periodic boundary conditions with period $4 \pi$ in all three 
directions.  The initial condition is the same as the one used by Kerr 
(see Section III of \cite{Kerr93}, and also \cite{HL06} for corrections of 
some typos in the description of the initial condition in \cite{Kerr93}). 
Following \cite{Kerr93}, we call the $x$-$y$ plane 
as the ``dividing plane'' and the $x$-$z$ plane as the ``symmetry plane''. 
There is one vortex tube above and below the dividing plane respectively. 
The term ``antiparallel'' refers to the anti-symmetry of the vorticity 
with respect to the dividing plane in the following sense: 
$\vec{\omega}(x,y,z) = -\vec{ \omega}(x,y,-z)$. Moreover, with respect to 
the symmetry plane, the vorticity is symmetric in its $y$ component and 
anti-symmetric in its $x$ and $z$ components. Thus we have 
$\omega_x(x,y,z)=-\omega_x(x,-y,z)$, $\omega_y(x,y,z)=\omega_y(x,-y,z)$ and
$\omega_z(x,y,z)=-\omega_z(x,-y,z)$. Here $\omega_x, \; \omega_y , \; \omega_z$
are the $x$, $y$, and $z$ components of vorticity respectively. These symmetries 
allow us to compute only one quarter of the whole periodic cell.

\begin{figure}
\begin{center} 
\includegraphics[width=8cm]{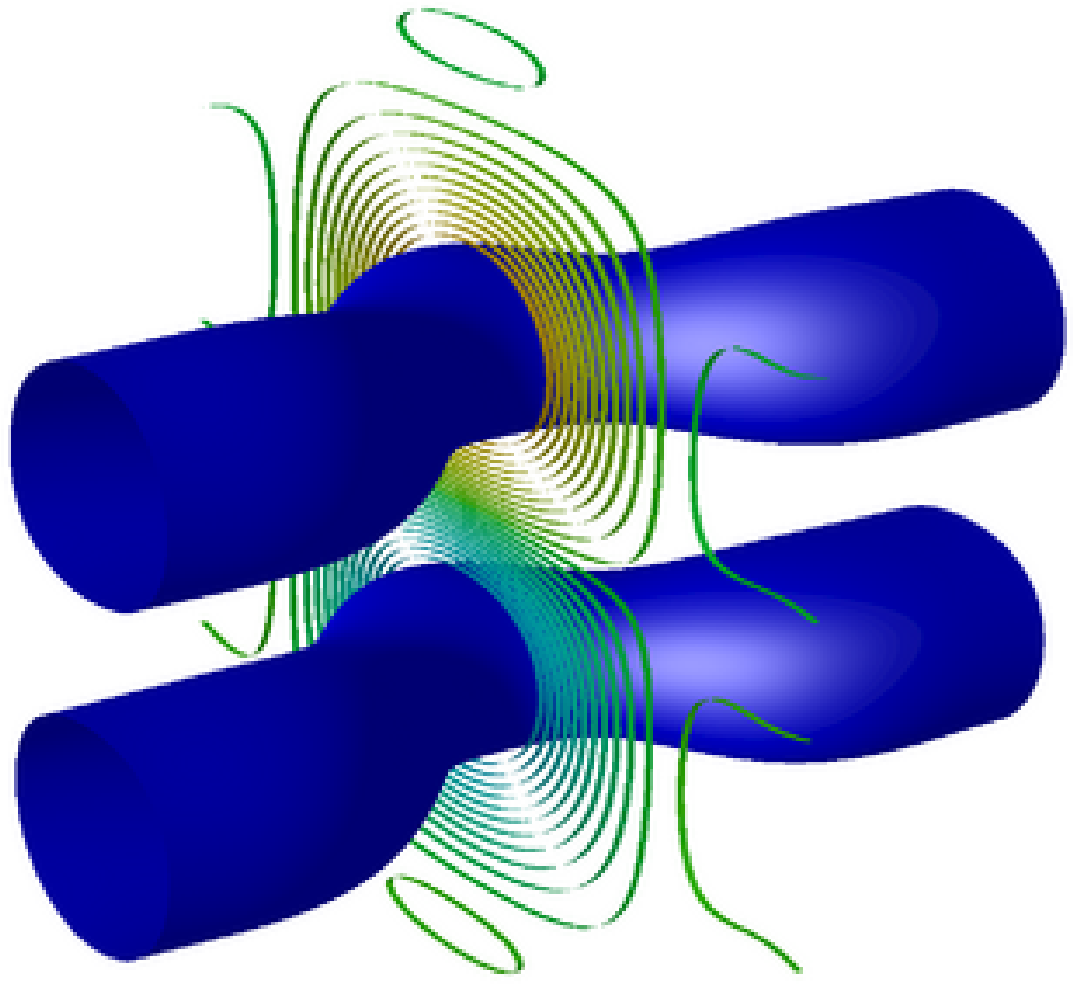}
\includegraphics[width=8cm]{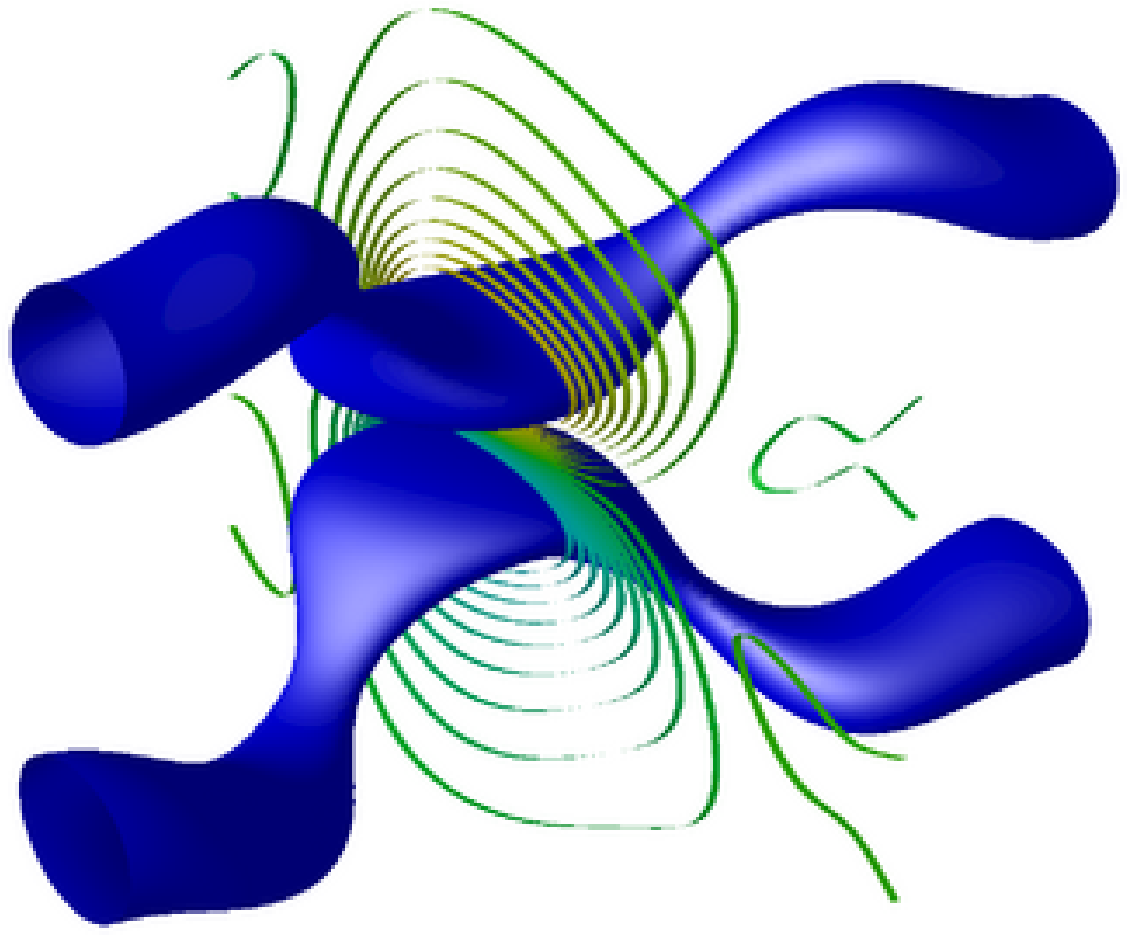}
\end{center}
\caption{The 3D view of the vortex tube for $t=0$ and $t=6$. The 
tube is the isosurface at $60\%$ of the maximum vorticity.
The ribbons on the symmetry plane are the contours
at other different values.
\label{fig.vorttube}}
\end{figure}

\begin{figure}
\begin{center}
\includegraphics[width=10cm]{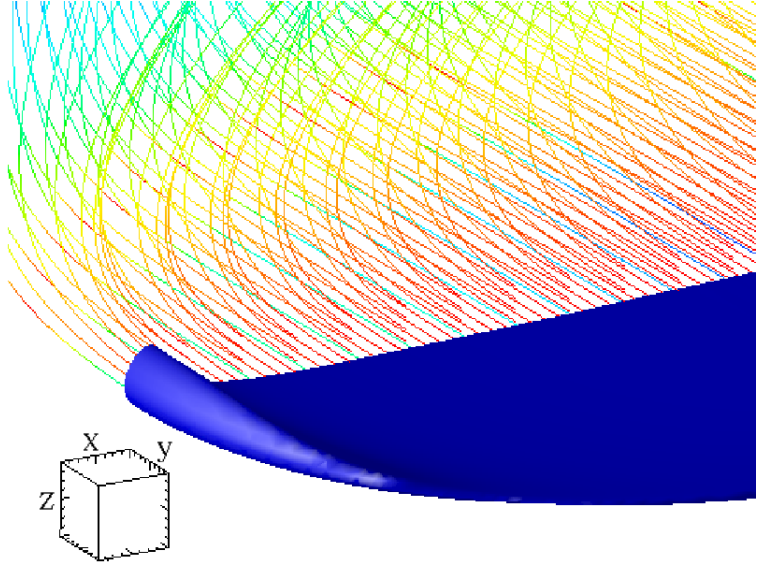}
\end{center}
\caption{The local 3D vortex structure and vortex lines around the maximum
vorticity at $t=17$. The size of the box on the left is
$0.075^3$ to demonstrate the scale of the picture.
\label{fig.local_struc_3d_17}}
\end{figure}

To compare the performance of the two pseudo-spectral methods, we will 
perform a careful convergence study for the two methods. To get a better 
idea how the solution evolves dynamically, we present the 3D plot of the 
vortex tubes at $t=0$ 
and $t=6$ respectively in Figure \ref{fig.vorttube}. As we can see, the 
two initial vortex tubes are very smooth and essentially symmetric. Due 
to the mutual attraction of the two antiparallel vortex tubes, the 
two vortex tubes approach to each one and experience severe deformation 
dynamically. By time $t=6$, there is already a significant flattening 
near the center of the tubes. In Figure \ref{fig.local_struc_3d_17}, we 
plot the local 3D vortex structure of the upper vortex tube at $t=17$.
By this time, the 3D vortex tube has essentially turned into a
thin vortex sheet with rapidly decreasing thickness. The vortex sheet
rolls up near the left edge of the sheet. It is interesting to note 
that the maximum vorticity is actually located near the rolled-up 
region of the vortex sheet.  

\subsection{Convergence study of the two pseudo-spectral methods in the
spectral space}

In this subsection, we perform a convergence study for the
two numerical methods using a sequence of resolutions.
For the Fourier smoothing method, we use the resolutions 
$768\times 512\times 1536$,
$1024\times 768\times 2048$, and $1536\times 1024\times 3072$
respectively. Except for the computation on the largest resolution
$1536\times 1024\times 3072$, all computations are carried out from 
$t=0$ to $t=19$. The computation on the final resolution
$1536\times 1024\times 3072$ is started from $t=10$ with the
initial condition given by the computation with the resolution
$1024\times 768\times 2048$. For the 2/3 dealiasing method, we use
the resolutions $512 \times 384 \times 1024$,
$768\times 512\times 1536$ and $1024\times 768\times 2048$
respectively. The computations using these three resolutions
are all carried out from $t=0$ to $t=19$. The time integration is 
performed using the classical fourth order Runge-Kutta method. 
Adaptive time stepping is used to satisfy 
the CFL stability condition with CFL number equal to $\pi/4$.

\begin{figure}
\begin{center}
\includegraphics[width=12cm,height=6cm]{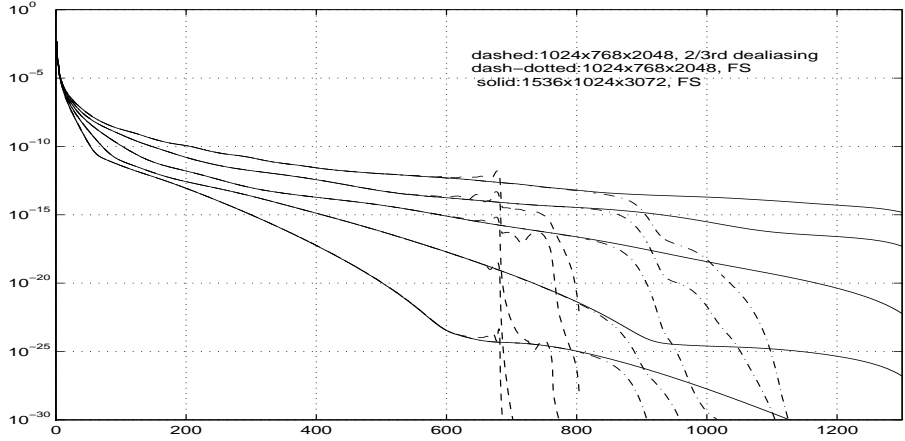}
\end{center}
\caption{The energy spectra versus wave numbers. We compare the energy
spectra obtained using the 
Fourier smoothing method with those using the 2/3 dealiasing method. The dashed
lines and the dashed-dotted lines are the energy spectra with the
resolution $1024\times 768\times 2048$ using the 2/3 dealiasing method and the
Fourier smoothing method, respectively. The solid lines are the energy spectra
obtained by the Fourier smoothing method with the highest resolution 
$1536\times 1024 \times 3072$. The times for the spectra lines
are at $t=15, 16, 17, 18, 19$ respectively.
\label{fig.energy-spec-comp}}
\end{figure}

\begin{figure}
\begin{center}
\includegraphics[width=12cm,height=6cm]{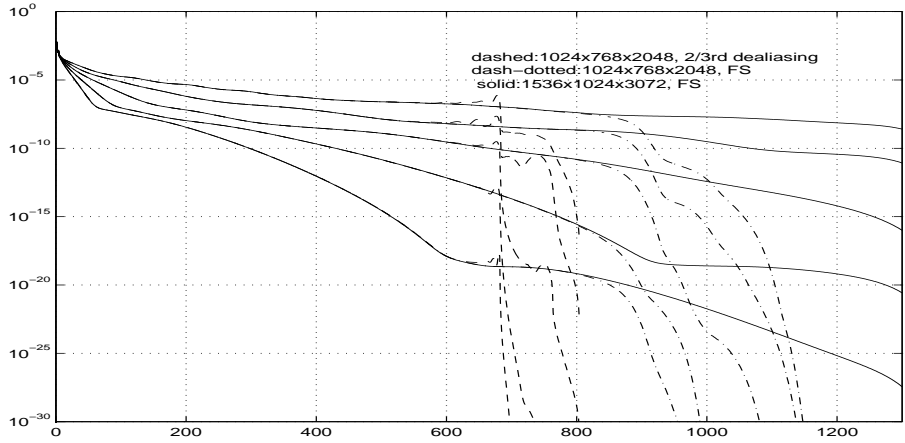}
\end{center}
\caption{The enstrophy spectra versus wave numbers. We compare the enstrophy 
spectra obtained using the 
Fourier smoothing method with those using the 2/3 dealiasing method. The dashed 
lines and dashed-dotted lines are the enstrophy spectra with the
resolution $1024\times 768\times 2048$ using the 2/3 dealiasing method and the
 Fourier smoothing, respectively. The solid lines are the enstrophy spectra
with resolution $1536\times 1024\times 3072$ obtained using the 
Fourier smoothing. The times for the spectra lines 
are at $t= 15, 16, 17, 18, 19$ respectively. 
\label{fig.enstrophy-spec-comp}}
\end{figure}

\begin{figure} 
\begin{center}
\includegraphics[width=12cm,height=6cm]{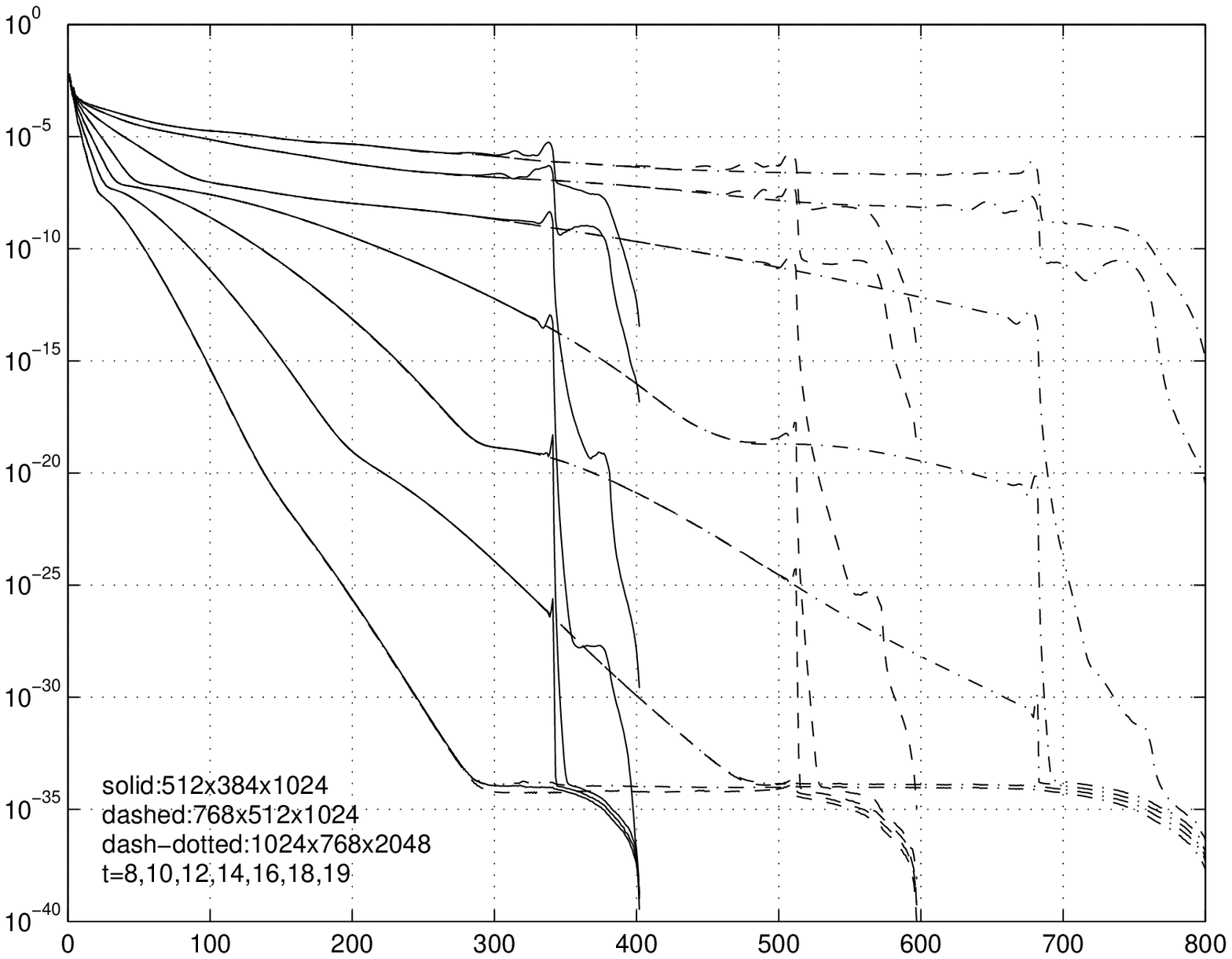}
\end{center}
\caption{Convergence study for enstrophy spectra obtained by
the 2/3 dealiasing method using different resolutions. The solid line 
is computed with
resolution $512\times 384\times 1024$, the dashed line is 
computed with resolution $786\times 512\times 1536$, and
the dashed-dotted line is computed with resolution 
$1024\times 768\times 2048$. The times for the lines 
from bottom to top are $t=8, 10, 12, 14, 16, 18, 19$.
\label{fig.enstrophy_spec_2}}
\end{figure}

In Figure \ref{fig.energy-spec-comp}, we compare the Fourier 
spectra of the energy obtained by using the $2/3$ dealiasing method
with those obtained by the Fourier smoothing method. 
For a fixed resolution $1024\times 768\times 2048$, we can see
that the Fourier spectra obtained by the Fourier smoothing
method retains more effective Fourier modes than those obtained 
by the $2/3$ dealiasing method. This can be seen by comparing the 
results with the corresponding computations using a higher 
resolution $1536\times 1024 \times 3072$. Moreover, the 
Fourier smoothing method does not give the spurious oscillations 
in the Fourier spectra which are present in the computations using 
the $2/3$ dealiasing method near the $2/3$ cut-off point. Similar 
convergence study has been made in the enstrophy spectra computed 
by the two methods. The results are given in Figures 
\ref{fig.enstrophy-spec-comp} and \ref{fig.enstrophy_spec_2}. 
They give essentially the same results.

\begin{figure}
\begin{center}
\includegraphics[width=12cm,height=6cm]{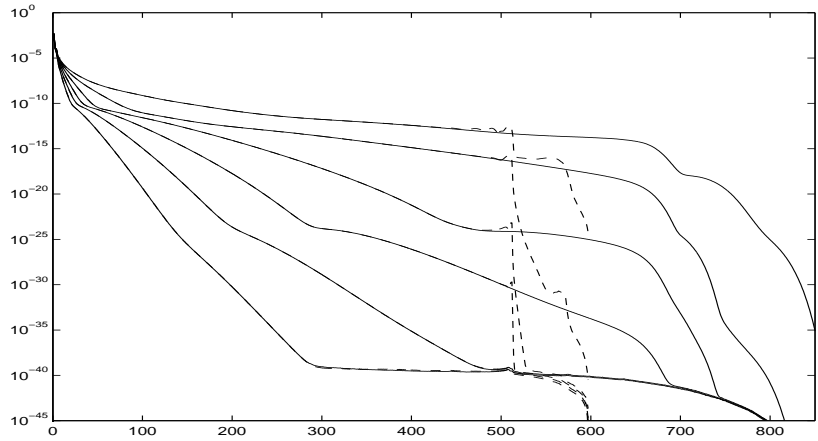}
\end{center}
\caption{The energy spectra versus wave numbers. We compare the energy
spectra obtained using the 
Fourier smoothing method with those using the 2/3 dealiasing method. The dashed
lines and solid lines are the energy spectra with the
resolution $768\times 512\times 1536$ using the 2/3 dealiasing method and the
 Fourier smoothing, respectively. The times for the spectra lines
are at $t=8, 10, 12, 14, 16, 18$ respectively.
\label{fig.energy-spec-comp-early}}
\end{figure}

We perform further comparison of the two methods using the
same resolution. In  Figure \ref{fig.energy-spec-comp-early},
we plot the energy spectra computed by the two methods using 
resolution $768\times 512\times 1536$. We can see that 
there is almost no difference in the Fourier spectra 
generated by the two methods in early times, $t=8, 10$,
when the solution is still relatively smooth. The difference  
begins to show near the cut-off point when the Fourier spectra 
rise above the round-off error level starting from $t=12$. 
We can see that the spectra computed by the 2/3 dealiasing 
method introduces noticeable oscillations near the 2/3 
cut-off point. The spectra computed by the Fourier
smoothing method, on the other hand, extend smoothly
beyond the 2/3 cut-off point. As we see from 
Figures \ref{fig.energy-spec-comp-early} and
\ref{fig.energy-spec-comp}, a significant portion 
of those Fourier modes beyond the 2/3 cut-off position 
are still accurate. This portion of the Fourier modes that
go beyond the 2/3 cut-off point is about $12\sim 15\%$ of total 
number of modes in each dimension. For 3D problems, the total 
number of effective modes in the Fourier smoothing method is 
about 20\% more than that in the 2/3 dealiasing method.
This is a very significant increase in the resolution for
a large scale computation. In our largest resolution, 
the effective Fourier modes in our Fourier smoothing method
are more than 320 millions, which has 140 millions more
effective modes than the corresponding 2/3 dealiasing method. 

\subsection{Comparison of the two methods in the physical space}

\begin{figure}
\begin{center}
\includegraphics[width=10cm,height=6cm]{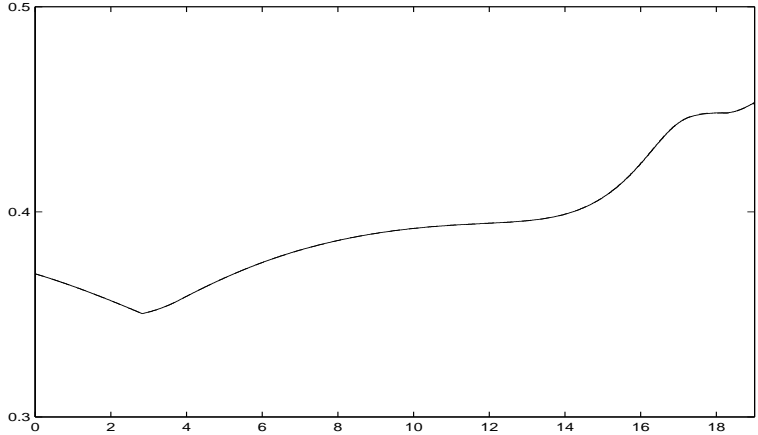}
\end{center}
\caption{Comparison of maximum velocity as a function of time
computed by two methods. The solid line represents the solution
obtained by the Fourier smoothing method,
and the dashed line represents the solution obtained by
the 2/3 dealiasing method.
The resolution is $1024\times 768\times 2048$ for both methods.
\label{fig.max-u-comp-1024}}
\end{figure}

\begin{figure}
\begin{center}
\includegraphics[width=10cm]{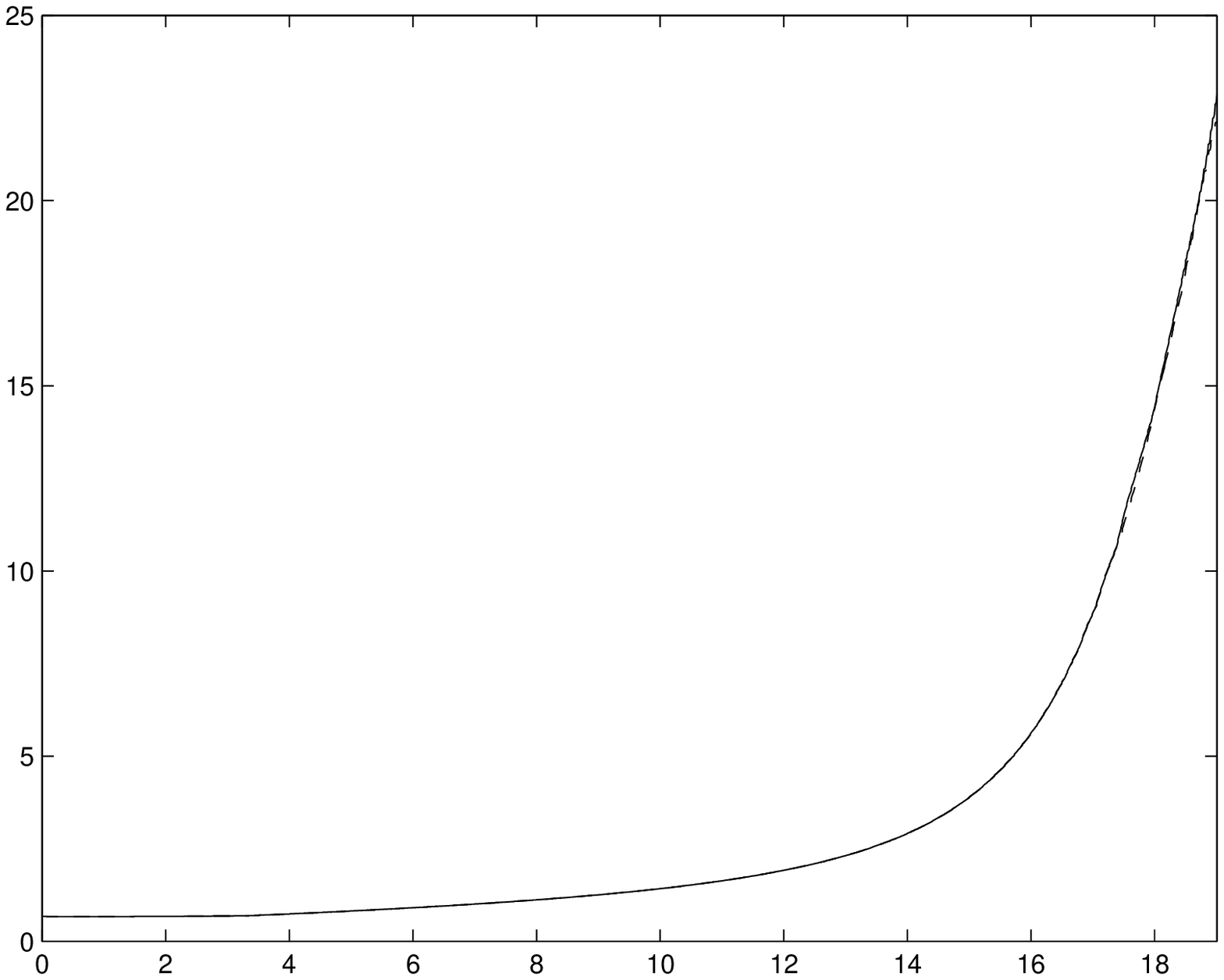}
\end{center}
\caption{Comparison of maximum vorticity as a function of time
computed by two methods. The solid line represents the solution
obtained by the Fourier smoothing method,
and the dashed line represents the solution obtained by
the 2/3 dealiasing method.
The resolution is $1024\times 768\times 2048$ for both methods.
\label{fig.max-vort-comp-1024}}
\end{figure}

\begin{figure}
\begin{center}
\includegraphics[width=10cm]{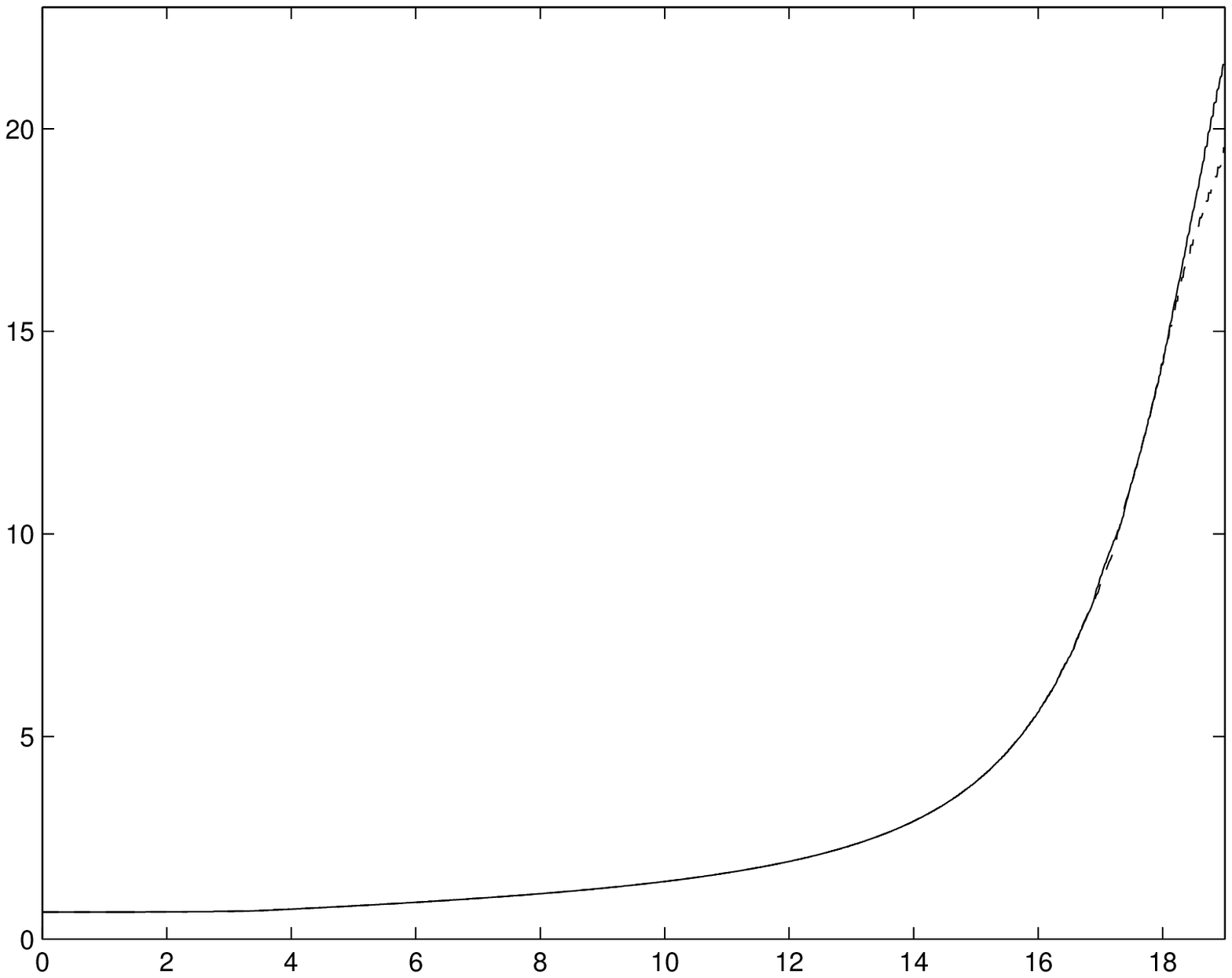}
\end{center}
\caption{Comparison of maximum vorticity as a function of time
computed by two methods. The solid line represents the solution
obtained by the Fourier smoothing method,
and the dashed line represents the solution obtained by
the 2/3 dealiasing method.
The resolution is $768\times 512\times 1024$ for both methods.
\label{fig.max-vort-comp-768}}
\end{figure}

\begin{figure}
\begin{center}
\includegraphics[width=9cm]{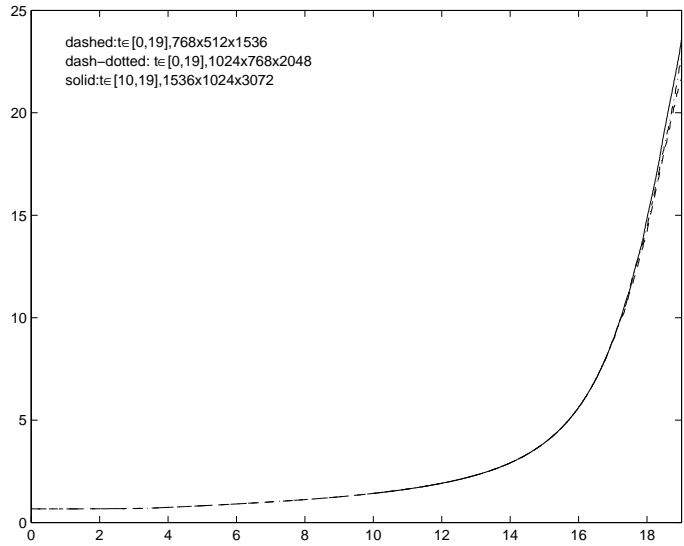}
\end{center} 
\caption{The maximum vorticity $\|\vec{\omega}\|_\infty$ in time 
computed by the Fourier smoothing method using
different resolutions.
\label{fig.omega}}
\end{figure}

\begin{figure}
\begin{center}
\includegraphics[width=9cm]{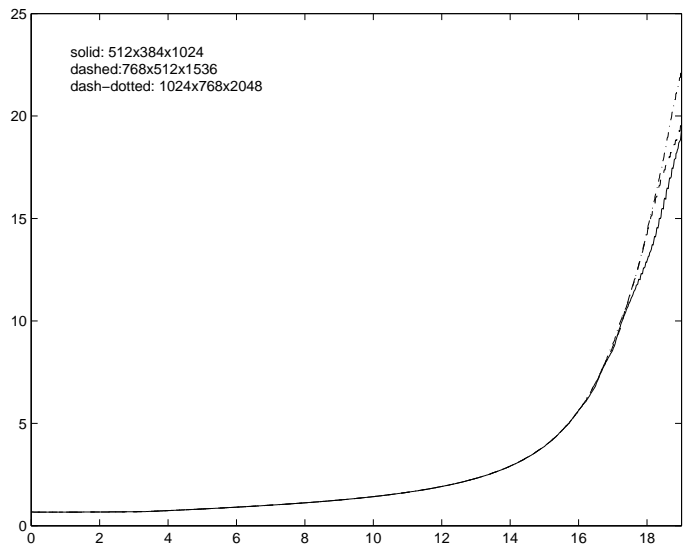}
\end{center}
\caption{The maximum vorticity $\|\vec{\omega}\|_\infty$ in time
computed by the 2/3 dealiasing method using different resolutions.
\label{fig.omega_2}} 
\end{figure}

Next, we compare the solutions obtained by the two methods 
in the physical space for the velocity field and the vorticity. 
In Figure \ref{fig.max-u-comp-1024}, 
we compare the maximum velocity as a function of time computed by 
the two methods using resolution $1024\times 768\times 2048$. The 
two solutions are almost indistinguishable. In Figure 
\ref{fig.max-vort-comp-1024}, we plot the maximum vorticity
as a function of time. The two solutions also agree reasonably 
well. However, the comparison of the solutions obtained by
the two methods at resolutions lower than $1024\times 768\times 2048$ 
shows more significant differences of the two methods, see
Figures \ref{fig.max-vort-comp-768}, \ref{fig.omega} and
\ref{fig.omega_2}.

To understand better how the two methods differ in their 
performance, we examine the contour plots of the axial vorticity 
in Figures \ref{fig.vort-cont-comp-1024-17},
\ref{fig.vort-cont-comp-1024-18} and \ref{fig.vort-cont-comp-1024-19}. 
As we can see, the vorticity computed by the 2/3 dealiasing method
already develops small oscillations at $t=17$. The oscillations grow 
bigger by $t=18$ (see Figure \ref{fig.vort-cont-comp-1024-18}), 
and bigger still at $t=19$ (see Figure \ref{fig.vort-cont-comp-1024-19}). 
We note that the oscillations in the axial vorticity contours concentrate
near the region where the magnitude of vorticity is close to zero. Thus
they have less an effect on the maximum vorticity. On the other hand, 
the solution computed by the Fourier smoothing method is still relatively 
smooth.

\begin{figure}
\begin{center}
\includegraphics[width=12cm,height=4cm]{pics/vort_cont_da_1024_17.epsf}
\includegraphics[width=12cm,height=4cm]{pics/vort_cont_fs36_1024_17.epsf}
\end{center}
\caption{Comparison of axial vorticity contours at $t=17$ 
computed by two methods. The upper picture is the solution
obtained by the 2/3 dealiasing method,
and the picture on the bottom is the solution obtained by
the Fourier smoothing method.
The resolution is $1024\times 768\times 2048$ for both methods.
\label{fig.vort-cont-comp-1024-17}}
\end{figure}
                                                                                          
\begin{figure}
\begin{center}
\includegraphics[width=12cm,height=4cm]{pics/vort_cont_da_1024_18.epsf}
\includegraphics[width=12cm,height=4cm]{pics/vort_cont_fs36_1024_18.epsf}
\end{center}
\caption{Comparison of axial vorticity contours at $t=18$ 
computed by two methods. The upper picture is the solution
obtained by the 2/3 dealiasing method,
and the picture on the bottom is the solution obtained by
the Fourier smoothing method.
The resolution is $1024\times 768\times 2048$ for both methods.
\label{fig.vort-cont-comp-1024-18}}
\end{figure}
                                                                                          
\begin{figure}
\begin{center}
\includegraphics[width=12cm,height=4cm]{pics/vort_cont_da_1024_19.epsf}
\includegraphics[width=12cm,height=4cm]{pics/vort_cont_fs36_1024_19.epsf}
\end{center}
\caption{Comparison of axial vorticity contours at $t=19$ 
computed by two methods. The upper picture is the solution
obtained by the 2/3 dealiasing method,
and the picture on the bottom is the solution obtained by
the Fourier smoothing method.
The resolution is $1024\times 768\times 2048$ for both methods.
\label{fig.vort-cont-comp-1024-19}}
\end{figure}
                                                                                          
%
To further demonstrate the accuracy of our computations we compare the 
maximum vorticity obtained by the Fourier smoothing method for three 
different resolutions: $768\times 512\times 1536$, 
$1024\times 768\times 2048$, and $1536\times 1024\times 3072$ respectively.
The result is plotted in Figure \ref{fig.omega}.
We have performed a similar convergence study for the
2/3 dealiasing method for the maximum vorticity. The result is
given in Figure \ref{fig.omega_2}. Two conclusions
can be made from this resolution study. First, by comparing
Figure \ref{fig.omega} with Figure \ref{fig.omega_2},
we can see that the Fourier smoothing method is indeed more accurate 
than the 2/3 dealiasing method for a given resolution. The 2/3
dealiasing method gives a slower growth rate in the maximum vorticity
with resolution $768\times 512\times 1536$. Secondly, the 
resolution $768\times 512\times 1536$ is not good enough to resolve 
the nearly singular solution at later times. On the other hand, we 
observe that the difference between the numerical solution obtained by
the resolution $1024\times 768\times 2048$ and that obtained by
the resolution $1536\times 1024\times 3072$ is relatively small.
This indicates that the vorticity is reasonably well-resolved
by our largest resolution $1536\times 1024\times 3072$. 

We have also performed a similar resolution study for the maximum
velocity in Figure \ref{fig.velocity}. The solutions obtained by
the two largest resolutions are almost indistinguishable,
which suggests that the velocity is well-resolved by
our largest resolution $1536\times 1024\times 3072$.

The resolution study given by Figures 
\ref{fig.max-vort-comp-768} and \ref{fig.omega_2} 
also suggests that the the computation obtained by the 
pseudo-spectral method with the 2/3 dealiasing rule using 
resolution $768\times 512\times 1536$ is significantly
under-resolved after $t=18$. It is interesting 
to note from Figure \ref{fig.max-vort-comp-768} that 
the computational results obtained by the two methods 
with resolution $768\times 512\times 1536$
begin to deviate from each other precisely around $t=18$. By 
comparing the result from Figure \ref{fig.max-vort-comp-768} 
with that from Figure \ref{fig.omega}, we confirm again 
that for a given resolution, the Fourier smoothing method
gives a more accurate approximation
than the 2/3 dealiasing method.

We remark that our numerical computations for the 3D incompressible
Euler equations using the Fourier smoothing method are very 
stable and robust. No high frequency instability has been observed 
throughout the computations. The resolution study we perform here is
completely based on the consideration of accuracy, not on stability.
This again confirms that the Fourier smoothing method offers a very 
stable and accurate computational method for the 3D incompressible flow.

\subsection{Does a finite time singularity develop?}

Before we conclude this section, we would like to have 
further discussions how to interpret the numerical results
we have obtained. Specifically, given the fast growth of
maximum vorticity, does a finite time singularity develop
for this initial condition? 

In \cite{Kerr93}, Kerr presented numerical evidence which suggested
a finite time singularity of the 3D Euler equations for the same
initial condition that we use in this paper. Kerr used a 
pseudo-spectral discretization with the 2/3 dealiasing rule
in the $x$ and $y$ directions, and a Chebyshev method in the $z$ 
direction with resolution of order $512\times 256 \times 192$. 
His computations showed that the growth of the peak vorticity, 
the peak axial strain, and the enstrophy production obey 
$(T-t)^{-1}$ with $T = 18.9$. In his recent paper \cite{Kerr05},
Kerr applied a high wave number filter to the data obtained 
in his original computations to ``remove the noise that 
masked the structures in earlier graphics'' presented in 
\cite{Kerr93}. With this filtered solution, he presented 
some scaling analysis of the numerical solutions up to 
$t=17.5$. Two new properties were presented in this recent
paper \cite{Kerr05}. First, the velocity field was shown to blow 
up like $O(T-t)^{-1/2}$ with $T$ being revised to $T=18.7$.
Secondly, he showed that the blowup is characterized by two 
anisotropic length scales, $\rho \approx (T-t) $ and 
$R \approx (T-t)^{1/2}$.

From the resolution study we present in Figure \ref{fig.omega},
we find that the maximum vorticity increases rapidly from the
initial value of $0.669$ to $23.46$ at the final time $t=19$,
a factor of 35 increase from its initial value. Kerr's
computations predicted a finite time singularity at $T=18.7$.
Our computations show no sign of finite time blowup of the 3D Euler
equations up to $T=19$, beyond the singularity time predicted by
Kerr. From Figures \ref{fig.vort-cont-comp-1024-17},
\ref{fig.vort-cont-comp-1024-18} and \ref{fig.vort-cont-comp-1024-19},
we can see that a thin layer (or a vortex sheet) is formed dynamically.
Beyond $t=17$, the vortex sheet has rolled up and traveled backward 
for some distance. With only 192 grid points along the $z$-direction,
Kerr's computations did not have enough grid points to resolve the
nearly singular vortex sheet that travels backward and away from
the $z$-axis. In comparison, we have 3072 grid points along the
$z$-direction. This gives about 16 grid points across the nearly 
singular layered structure at $t=18$ and about $8$ grid points at 
$t=19$.

\begin{figure}
\begin{center}
\includegraphics[width=8cm]{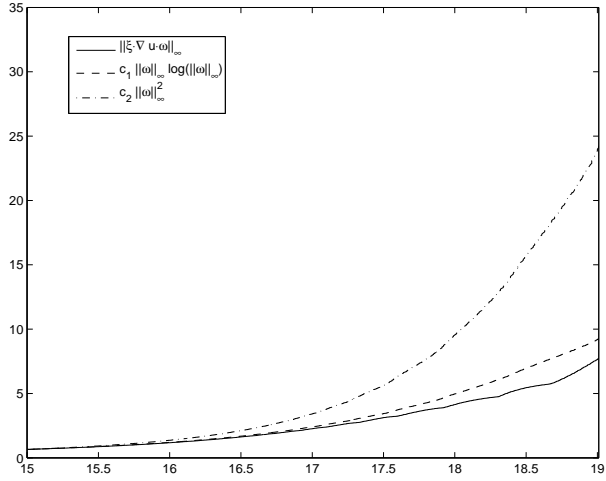}
\end{center}
\caption{Study of the vortex stretching term in time. This computation is 
performed by the Fourier smoothing method with resolution $1536\times
1024\times 3072$. We take $c_1 = 1/8.128$, $c_2 = 1/23.24$ to match the
same starting value for all three plots.
\label{fig.growth_rate}}
\end{figure}
                                                                                
\begin{figure}
\begin{center}
\includegraphics[width=8cm]{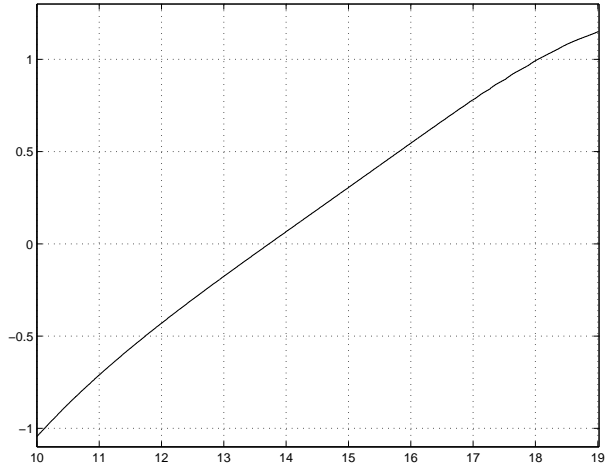}
\end{center}
\caption{The plot of $ \log \log \|\omega\|_\infty$ vs time. 
This computation is performed by the Fourier smoothing 
method with resolution $1536\times 1024 \times 3072$.
\label{fig.omega_loglog}}
\end{figure}
                                                                                
In order to understand the nature of the dynamic growth in
vorticity, we examine the degree of nonlinearity in the vortex 
stretching term. In Figure \ref{fig.growth_rate}, we plot
the quantity,
$\|\xi \cdot \nabla \vec{u} \cdot \vec{\omega}\|_\infty$,
as a function of time, where $\xi$ is the unit vorticity vector.
If the maximum vorticity indeed blew up like $O((T-t)^{-1})$,
as alleged in \cite{Kerr93}, this quantity should have been
quadratic as a function of maximum vorticity. We find that
there is tremendous cancellation in this vortex stretching
term. It actually grows slower than
$C\|\vec{\omega} \|_\infty \log (\|\vec{\omega} \|_\infty )$,
see Figure \ref{fig.growth_rate}. It is easy to show that 
such weak nonlinearity in vortex stretching would imply 
only doubly exponential growth in the maximum vorticity.
Indeed, as demonstrated by Figure \ref{fig.omega_loglog}, the
maximum vorticity does not grow faster than doubly exponential
in time. In fact, a closer inspection reveals that the location 
of the maximum vorticity has moved away from the dividing plane
for $t \geq 17.5$. This implies that the compression mechanism
between the two vortex tubes becomes weaker toward the end of the
computation, leading to a slower growth rate in maximum vorticity
\cite{HL06}.
                                                                                
\begin{figure}
\begin{center}
\includegraphics[width=9cm]{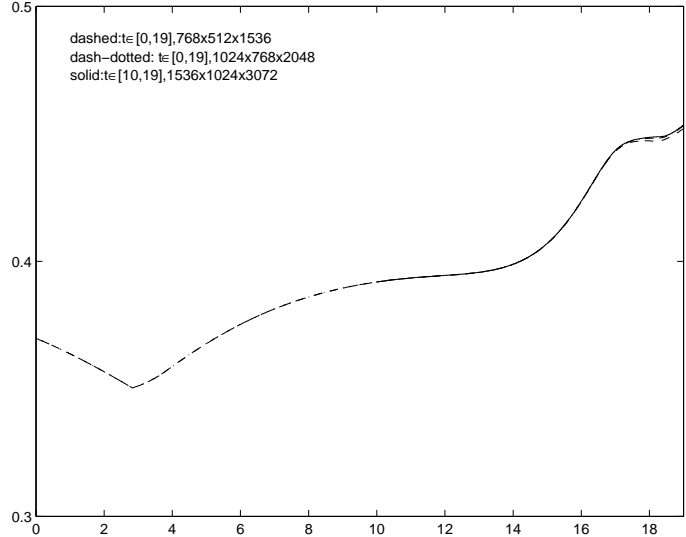}
\end{center}
\caption{Maximum velocity $\|\vec{u}\|_\infty$ in time
computed by the Fourier smoothing method using different resolutions.
\label{fig.velocity}} 
\end{figure} 

Another important evidence which supports the non-blowup of
the solution up to $t=19$ is that the maximum velocity remains
bounded, see Figure \ref{fig.velocity}. This is in contrast
with the claim in \cite{Kerr05} that the maximum velocity 
blows up like $O(T-t)^{-1/2}$ with $T=18.7$. With the velocity
field being bounded, the local non-blowup criteria of Deng-Hou-Yu
\cite{DHY05a,DHY05b} can be applied, which implies that the solution
of the 3D Euler equations remains smooth at least up to $T=19$,
see also \cite{HL06}.

\section{Conclusion Remarks}

In this paper, we have performed a systematic convergence study of
the two pseudo-spectral methods. The first pseudo-spectral method
uses the traditional 2/3 dealiasing rule, while the second
pseudo-spectral method uses a high order Fourier smoothing.
The Fourier smoothing method is designed to cut
off the high frequency modes smoothly while retaining a significant
portion of the Fourier modes beyond the 2/3 cut-off point in
the spectral space. We apply both methods to compute nearly
singular solutions of the 1D Burgers equation and the 3D 
incompressible Euler equations. In the case of the 1D Burgers
equation, we can obtain a very accurate approximation of the
exact solution sufficiently close to the singularity time
with 13 digits of accuracy. This allows us to estimate
the numerical errors of the two methods accurately and provides
a solid ground in our convergence study. In our study of
the 3D incompressible Euler equations, we use the highest
resolution that we can afford to perform our convergence study.
In both cases, we demonstrate convincingly that the Fourier smoothing
method gives a more accurate approximation than the 2/3 dealiasing method.

Our extensive convergence studies in both physical and spectral
spaces show that the Fourier smoothing method offers
several advantages over the 2/3 dealiasing method when computing
a nearly singular solution. First of all, the error in the
Fourier smoothing method is highly localized near the region
where the solution is most singular. The error in the smooth region 
is several orders of magnitude smaller than that near the 
``singular'' region. The 2/3 dealiasing method, on the other hand, 
has a wide spread pointwise error distribution, and produces relatively 
large oscillations even in the smooth region. Secondly, for the same 
resolution, the Fourier smoothing method offers a 
more accurate approximation to the physical solution than the 
2/3 dealiasing method. Our numerical study shows that for a 
given resolution, the Fourier smoothing method retains 
about $12\sim 15\%$ more effective Fourier modes than the 2/3 
dealiasing method in each dimension. For a 3D problem, the gain 
is as large as 20\%. This gain is quite significant in a 
large scale computation. Thirdly, the Fourier smoothing 
method is very stable and robust when computing nearly singular 
solutions of fluid dynamics equations. Spectral convergence is 
clearly observed in all our computational experiments without 
suffering from the Gibbs phenomenon. Moreover, there is no additional 
computational cost in implementing the Fourier 
smoothing method. We have also implemented the Fourier
smoothing method for the incompressible 3D Navier-Stokes equations
and observed a similar performance.

We have applied both spectral methods to study the 
potentially singular solution of the 3D Euler equation using 
the same initial condition as Kerr \cite{Kerr93}. Both the
Fourier smoothing method and the 2/3 dealiasing
method give qualitatively the same result except that the
2/3 dealiasing method suffers from the Gibbs phenomenon
and produces relative large oscillations at late times.  
Our convergence study in both the physical and spectral spaces 
shows that the maximum vorticity does not grow faster
than double exponential in time and the maximum velocity field
remains bounded up to $T=19$, beyond the singularity time
$T = 18.7$ alleged in \cite{Kerr93,Kerr05}. Tremendous
cancellation seems to take place in the vortex stretching
term. The local geometric regularity of the vortex lines
near the region of the maximum vorticity seems to be 
responsible for this dynamic depletion of vortex stretching
\cite{CFM96,DHY05a,DHY05b,HL06}.

\vspace{0.2in}
\noindent
{\bf Acknowledgments.}
We would like to thank Prof. Lin-Bo Zhang from the Institute of
Computational Mathematics in Chinese Academy of Sciences (CAS) for 
providing us with the computing resource to perform this large
scale computational project. Additional computing resource was
provided by the Center of High Performance Computing in CAS. We
also thank Prof. Robert Kerr for providing us with his Fortran
subroutine that generates his initial data. This work was in part
supported by NSF under the NSF FRG grant DMS-0353838 and ITR
Grant ACI-0204932. Part of this work was done while Hou visited
the Academy of Systems and Mathematical Sciences of CAS in the
summer of 2005 as a member of the Oversea Outstanding Research
Team for Complex Systems. Li was supported by the National Basic 
Research Program of China under the grant 2005CB321701. Finally,
we would like to thank Professors Alfio Quarteroni, Jie Shen, and
Eitan Tadmor for their valuable comments on our draft manuscript.  

\bibliographystyle{amsplain}
\bibliography{bib}

\end{document}